\documentclass[11pt]{amsart}

\usepackage[T1]{fontenc}
\usepackage{lmodern,microtype}
\usepackage[margin=1.5in]{geometry}
\usepackage{amssymb,mathtools}
\usepackage{graphicx,tikz}
\usepackage[colorlinks=true,linkcolor=blue,citecolor=blue,urlcolor=blue]{hyperref}

\mathtoolsset{showonlyrefs}
\numberwithin{equation}{section}

\newtheorem{theorem}{Theorem}[section]
\newtheorem{corollary}[theorem]{Corollary}
\newtheorem{lemma}[theorem]{Lemma}
\newtheorem{remark}[theorem]{Remark}

\makeatletter
\newenvironment{prooflemma}[1]{%
  \def\thmhead##1##2##3{%
    {\hypersetup{linkcolor=black}%
     \hyperlink{#1}{\thmhead@plain{##1}{##2}{##3}}}}%
  \begin{lemma}%
}{\end{lemma}}
\makeatother

\newcommand{\abbr}[1]{{\sc\lowercase{#1}}}
\newcommand{\DLA}{\abbr{dla}}
\newcommand{\PP}{\mathbf P}
\newcommand{\EE}{\mathbf E}
\newcommand{\R}{\mathbb R}
\newcommand{\Z}{\mathbb Z}
\newcommand{\ind}{\mathbf 1}
\newcommand{\cF}{\mathcal F}
\newcommand{\dist}{\operatorname{dist}}

\title[\DLA{}: boundary, vacancies, and shape exclusion]
{Diffusion-limited aggregation:\\
boundary, vacancies, and shape exclusion}

\author{Amir Dembo}
\email{adembo@stanford.edu}
\author{Matthew Junge}
\email{Matthew.Junge@baruch.cuny.edu}

\hypersetup{
  pdftitle={Diffusion-limited aggregation: Boundary, vacancies and shape exclusion},
  pdfauthor={Amir Dembo and Matthew Junge}
}

\begin{document}
\begin{abstract}
For external Diffusion-Limited Aggregation (\DLA{}) on $\Z^d$ with $d \geq 2$, we prove that the number of boundary edges grows linearly.
Consequently, a positive fraction of the aggregate lies outside
every compact region rescaled to have volume equal to the
number of particles.
We also show that the probability that the limiting aggregate covers more than a $1-\alpha$ fraction of a cube decays exponentially in its volume for some $\alpha>0$.
For Euclidean \DLA{} in every dimension $d\ge2$, we show that
each particle has approximately two near neighbors on average.
We then deduce that a positive fraction of the planar Euclidean \DLA{} particles 
lie outside any compact region scaled to the optimal sphere packing volume.
\end{abstract}
\maketitle

\section{Introduction}

Diffusion-Limited Aggregation (\DLA{}) was introduced by Witten and
Sander as a model of fractal growth patterns seen in natural systems
in which particles diffuse from far away and attach to a growing aggregate \cite{Witten1981Diffusion}.
It has since been a central focus of the physics literature with thousands of publications
devoted to its study \cite{Sander2000Diffusion}.
In contrast, relatively few rigorous results are known (see the summary in Section~\ref{sec:past}).

\begin{figure}
  \centering
  \includegraphics[width=0.45\textwidth]{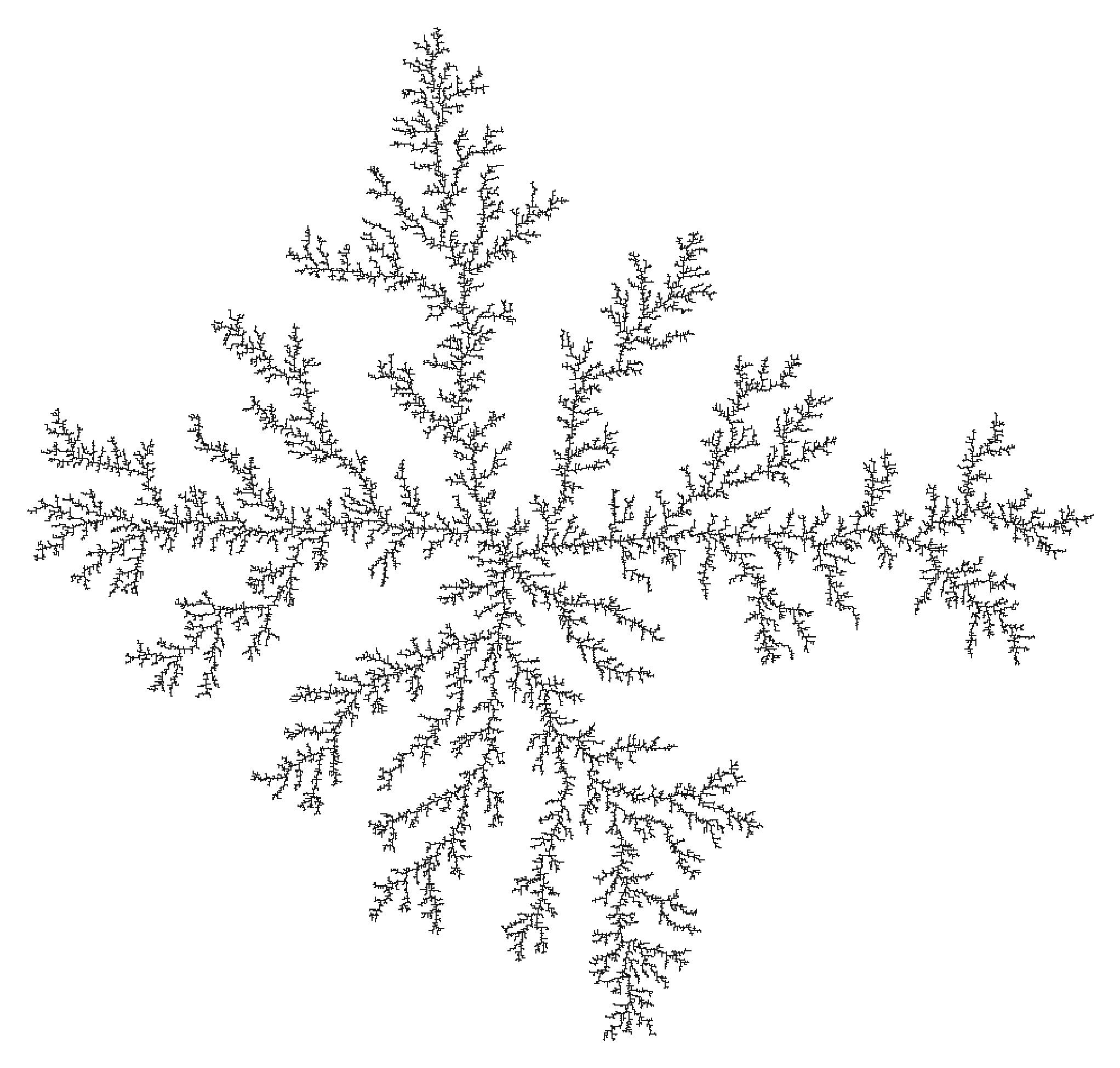} \qquad 
  \includegraphics[width=0.4\textwidth]{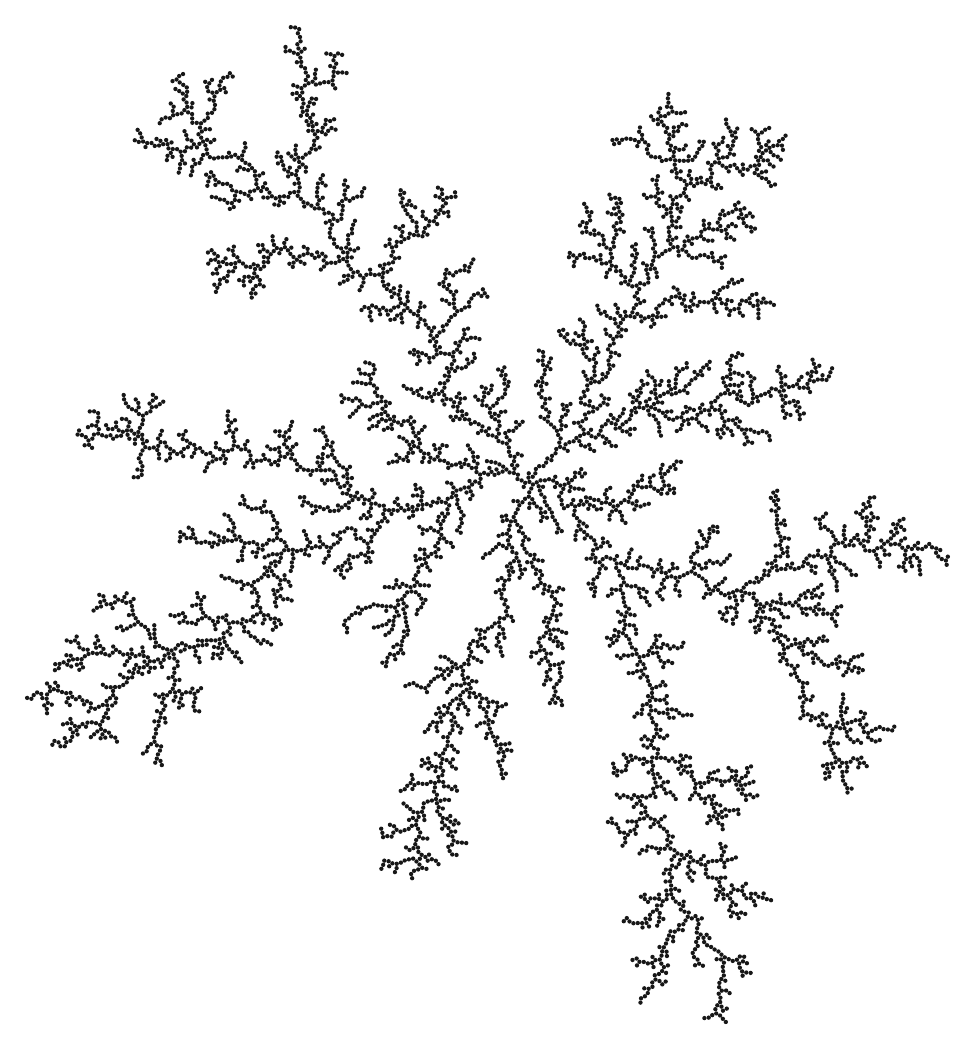}
  \caption{Simulations of (left) \DLA{} on $\mathbb{Z}^2$ with
  $n=50{,}000$ and (right) planar Euclidean \DLA{} with $n=5{,}000$.}
  \label{fig:dla-simulation}
\end{figure}

We consider two versions of this model. In lattice \DLA{} on $\mathbb Z^d$, the aggregate starts from $\{0\}$ and each incoming particle performs a simple random walk 
until it reaches an unoccupied site adjacent to the aggregate. 
In Euclidean \DLA{} on $\mathbb R^d$, the aggregate starts
from a unit ball centered at the origin. Each incoming unit ball has
its center follow Brownian motion until the ball touches the aggregate. 
See Figure~\ref{fig:dla-simulation} for a simulation of the aggregates of both models.

Theorems~\ref{thm:boundary}, \ref{thm:vacancies}, and~\ref{thm:attachment}
are our main probabilistic contributions.
Theorem~\ref{thm:boundary} proves that the edge boundary of lattice
\DLA{} grows linearly. Theorem~\ref{thm:vacancies} bounds the probability
that a large lattice cube ever becomes nearly full and gives uniform
control of permanent vacancies. Theorem~\ref{thm:attachment} gives uniform bounds on near contacts
in Euclidean \DLA{} and shows that the graph joining sufficiently
close balls has average degree arbitrarily close to $2$.
These results hold in every dimension $d\ge2$.
Deterministic boundary comparisons give the lattice shape exclusion in
Corollary~\ref{cor:shapes}. Sphere-packing comparisons give an analogous
shape exclusion at optimal density for planar Euclidean \DLA{}, 
see Corollary~\ref{cor:euclidean}.

\subsection{Lattice \DLA{}}\label{subsec:lattice-dla}

We first consider \emph{external site \DLA{}} on $\Z^d$,
$d\ge2$. 
For a finite nonempty connected set $A\subset\Z^d$, define
\[
 \partial A\coloneqq \{y\notin A\mid |y-x|_1=1\text{ for some }x\in A\},
 \qquad \overline A\coloneqq A\cup\partial A.
\]
For simple random walk $(S_j)$ started outside $\overline A$, consider the first entrance time 
\[
 \tau_A\coloneqq \inf\{j\ge0\mid S_j\in\overline A\}
\]
and the associated {harmonic measure from infinity}, namely the 
distribution
\begin{equation}\label{def:lattice-harmonic-measure}
 \mu_A(y)\coloneqq \lim_{|z|\to\infty}
 \PP_z(S_{\tau_A}=y\mid\tau_A<\infty),\qquad y\in\partial A.
\end{equation}
For existence of the limit when $d=2$ and $d \ge 3$, see for example 
\cite[Propositions~6.6.1 and~6.5.4]{Lawler2010Random}, respectively
(if $d=2$, the event $\tau_A<\infty$ holds with probability one). Starting the lattice \DLA{}
at $A_1 \coloneqq \{0\}$, sequentially in $n$, conditional on $A_n$, the next site 
$X_{n+1}$ is drawn from $\mu_{A_n}$ and added to the aggregate so that 
$A_{n+1}\coloneqq A_n\cup\{X_{n+1}\}$.

Our first result quantifies the irregularity of the \DLA{} aggregate through 
the growth rate of the \emph{edge boundary} 
\[
 \partial_e A\coloneqq \bigl\{\{x,y\}\mid x\in A,\ y\notin A,
                                      \ |x-y|_1=1\bigr\} .
\]
Note that $\partial_e A$ counts both edges to finite internal holes and those to the external infinite component of the complement. 

In what follows, we use the notation $|A|$ for cardinality or Lebesgue measure according to context.
The positive constants $\alpha$ appearing in different statements
need not have the same value. Noting that $|\partial_e A_n| \le 2 (d-1) n+2$ deterministically, our first result establishes a linear lower bound. 

\begin{theorem}\label{thm:boundary}
For every $d\ge2$, there is a constant $b_d>0$ such that almost surely
\[
 \liminf_{n\to\infty}
 \frac{|\partial_e A_n|}{n}\ge b_d.
\]
\end{theorem}

Bubeck in 2015 remarked that not even the ball had been ruled out as a limiting shape for \DLA \ 
\cite{Bubeck2015Some}. The linear lower bound in Theorem~\ref{thm:boundary} is much larger 
than the surface-order scale $n^{(d-1)/d}$ of a region with smooth boundary, allowing us to 
exclude any compact shape of unit volume. 
Specifically, for a bounded region $R\subset\R^d$ of unit volume, we define $R_n$ to be the lattice points contained in the rescaling of $R$ to have volume $n$:
\begin{equation}\label{def:Rn}
 R_n\coloneqq \Z^d\cap (n^{1/d} R). 
\end{equation}
So long as $|\partial R|=0$, i.e.\ $R$ is \emph{Jordan measurable}, we have $|R_n|=n+o(n)$. 
The next result also applies to bounded Jordan measurable
sets since taking their closure preserves unit volume.

\begin{samepage}
\begin{corollary}\label{cor:shapes}
For every $d\ge2$, there is a constant $\alpha=\alpha(d)>0$ such that almost surely, 
for every compact set $R\subset\R^d$ of unit volume
\begin{equation}\label{eq:exclude-shape}
 \liminf_{n\to\infty}\frac{|A_n\setminus R_n|}{n}\ge\alpha.
\end{equation}
\end{corollary}
\end{samepage}

Since the probability one event in Corollary~\ref{cor:shapes} is common
to all compact regions of unit volume, the set $R$ may be chosen
\emph{after observing the entire \DLA{} realization}. 
More generally, 
on that same probability one event, our proof shows that 
\begin{equation}\label{eq:exclude-small-bdry}
\liminf_{n\to\infty}\frac{|A_n\setminus B_n|}{n} \ge\alpha,
\end{equation}
for every sequence of finite lattice sets $B_n$ with $|B_n|=n+o(n)$ and 
$|\partial_e B_n|=o(n)$.
For example, $B_n$ could consist of unions of cubes with side length
$L_n\to\infty$. The selected cubes may form patterns that change with $n$ and are constructed with full
knowledge of the entire \DLA{} realization. Even so, a positive
fraction of the aggregate lies outside their union for all
sufficiently large $n$.

Another consequence of Corollary~\ref{cor:shapes} is that external \DLA{} differs 
from \emph{internal} \DLA{}, which is known to have the ball as its limiting shape \cite{Lawler1992Internal},
by a positive fraction of its volume.
To our knowledge, this is the first rigorous result establishing
a macroscopic distinction between the occupied sets of these two models.

A \emph{lattice cube} $B$ of side $m$ is a translate of
$\{1,\ldots,m\}^d$. Our next result concerns \emph{permanent vacancies}, i.e. the complement of $A_\infty\coloneqq\bigcup_{n\ge1}A_n$.
We prove that the probability $A_\infty$ misses at most a fixed
small fraction of the points in $B$ is exponentially small in $|B|$, 
complementing Eberz-Wagner's finding that planar \DLA{}
leaves infinitely many permanent holes
\cite{EberzWagner1999Discrete}.
We also prove a local version of Theorem~\ref{thm:boundary}.
If $A_\infty$ occupies a positive fraction of $B$, then the edge
boundary of $A_\infty\cap B$ has order $|B|$, except on an event
of exponentially small probability in $|B|$.

\begin{theorem}\label{thm:vacancies}
The following hold for every $d\ge2$.
\newline
(a). For some $\alpha=\alpha(d) >0$ and
every lattice cube $B$ of side at least $1/\alpha$,
\begin{equation}\label{eq:vacancy-tail}
 \PP\bigl(|B\setminus A_\infty|\le \alpha|B|\bigr)\le e^{-\alpha|B|}.
\end{equation}
(b). For some $\gamma=\gamma(d)>0$, any $\rho > 0$ 
and all lattice cubes $B$ of side at least $1/(\gamma \rho)$ 
\begin{equation}\label{dfn:loc-edge-bdry}
\PP\Bigl(
 |B\cap A_\infty|\ge\rho|B|,
 \quad
 |\partial_e(B\cap A_\infty)|\le \gamma |B \cap A_\infty|
 \Bigr)
 \le e^{- \gamma \rho^2|B|}.
\end{equation}
\end{theorem}
We note in passing that Theorem~\ref{thm:vacancies} applies for more general
collections of non-random finite sets, see Remark \ref{rmk:general-boxes}. 
Moreover, for $L:=\exp\{\alpha m^d/[2(d+1)]\}$, a union bound over 
the at most $(2L+1)^{d+1}$ lattice cubes in $[-L,L]^d$, followed by
Borel--Cantelli, shows that almost surely, for all sufficiently
large integers $m$,
\[
 |B\setminus A_\infty|\ge\alpha|B|
\]
simultaneously for every such cube $B$ of side at least $m$. Similarly, for 
$L=\exp\{\gamma \rho^2 m^d/[2(d+1)]\}$ and any $\rho>0$ we have that almost surely,
for all large integers $m$ and every lattice cube $B \subset [-L,L]^d$ 
of side at least $m$, 
\[
\frac{|\partial_e (A_\infty \cap B)|}{|A_\infty \cap B|} > \gamma
\quad {\rm or} \quad \frac{|A_\infty \cap B|}{|B|} < \rho \,.
\]
When $A_\infty\cap B=\varnothing$, only the second alternative is
needed, and it holds since $\rho>0$.

\subsection{Euclidean \DLA{}}
In \emph{Euclidean \DLA{}}, unit balls with disjoint
interiors attach according to Brownian harmonic measure from infinity.
To define the model, write $B(x,r)$ for the open Euclidean ball with center $x$ and radius $r$.
A set $C$ is called \emph{$2$-separated} if any two distinct points in $C$ are at 
distance at least $2$. With any such finite nonempty $2$-separated set $C$, we associate 
\begin{equation}\label{eq:exc-set}
 \mathcal O_C\coloneqq\bigcup_{a\in C}\overline B(a,2),
\end{equation}
and let $\mu_C$ denote the Brownian harmonic measure from infinity on $\partial\mathcal O_C$.
More precisely, consider the first entrance time 
$\tau_{C} \coloneqq\inf\{t\ge0\mid W(t) \in \mathcal O_C\}$
for Brownian motion $(W(t), t \ge 0)$ started at $z \notin \mathcal O_C$,
setting the total variation limit
\begin{equation}\label{def:harmonic-measure}
 \mu_C \coloneqq\lim_{|z|\to\infty}
 \mathcal L_z(W(\tau_{C}) \mid \tau_{C} <\infty).
\end{equation}
As before, there is no need to condition when $d=2$ and for the existence of this 
limit, see for example \cite[Proof of Theorem 3.46]{Morters2010Brownian}.

Growing the Euclidean \DLA{} at step $n$ one uses 
$\mu_{C_n}$ for the law of an incoming ball center at its first contact with the 
unit balls centered at $C_n$. That is, starting the 
aggregate at $C_1\coloneqq\{0\}$, sequentially in $n$, conditional on the $2$-separated set $C_n$ of centers of the $n$ unit balls, one samples $X_{n+1}$ from $\mu_{C_n}$ and sets 
$C_{n+1}\coloneqq C_n\cup\{X_{n+1}\}$. Note that each new unit ball touches the 
existing aggregate and has interior disjoint from all earlier balls. In this context, 
the attachment positions 
\begin{align}\label{dfn:good-attachment}
 E_\eta(C)\coloneqq 
 \bigcup_{\substack{a,b\in C\\a\ne b}}
 \{x \in \partial\mathcal O_C 
 \mid x \in \partial \overline B(a,2),\ |x-b|\le2+\eta\},
\end{align}
imply that the incoming ball
touches one existing ball and leaves a small boundary gap of at most
$\eta$ from another ball. More generally, the boundary gap between 
two unit balls of centers $x,y \in C$ is $|x-y|-2 \ge 0$.

Our next result shows that such \emph{almost-double contact}, where
$X_{k+1} \in E_\eta(C_k)$, is unlikely, uniformly over all stages of 
growth and with respect to the past realization. To fully state our result, 
we also define for $\eta>0$ and any finite set of centers $C$, the 
\emph{proximity graph} $G_\eta(C)$ on $C$ whose edges correspond 
to all $x \ne y \in C$ with $|x-y| \le 2 + \eta$. We also let $\mathcal{E}(G_\eta(C))$
denote the edge set. 

\begin{theorem}\label{thm:attachment}
For every $d\ge2$, there exists $c_d<\infty$ such that:

\noindent\textup{(a)} For any
finite nonempty $2$-separated set $C\subset\R^d$ and $0<\eta<1$,
\begin{equation}\label{eq:unif-attachment}
 \mu_C\bigl(E_\eta(C)\bigr)
 \le c_d\sqrt\eta.
\end{equation}
\noindent\textup{(b)} For each $0<\eta<1$, almost surely
\begin{equation}\label{eu-eq:degree-bound}
 \limsup_{n \to \infty} 
 \frac{| \mathcal{E}(G_\eta(C_n)) |}{n}
\le 1 + c_d \sqrt{\eta} .
\end{equation}
When $d\ge3$, the $c_d \, \sqrt{\eta}$ terms on the 
\abbr{rhs} of \eqref{eq:unif-attachment}-\eqref{eu-eq:degree-bound} 
improve to $c_d \, \eta$.
\end{theorem}

Part~\textup{(a)} shows that the probability of an almost-double contact can, uniformly over growth histories,
be made arbitrarily small. See Figure~\ref{fig:thin-attachment} for the geometry behind
this estimate.

Part~\textup{(b)} states that the number of proximity edges
is asymptotically at most $(1+O(\sqrt{\eta}))n$.
The proximity graph $G_\eta(C_n)$ which contains all connected contacts,
has at least $n-1$ edges, so the constant $1$ is optimal.
Equivalently, the average number of nearly touching neighbors
per ball is asymptotically at most $2+O(\sqrt{\eta})$.
Together with connectivity, this implies that a positive
fraction of the balls have at most two neighbors within a
fixed positive gap.
This geometric sparsity is analogous in spirit to
Theorem~\ref{thm:boundary}, although it does not by itself give
a lower bound on the external boundary of the exclusion region
$\mathcal O_{C_n}$.

\begin{figure}
\centering
\begin{minipage}[c]{.57\textwidth}
\centering
\begin{tikzpicture}[scale=.44]
  \definecolor{thinred}{RGB}{190,55,45}
  \def\gapwidth{0.18}
  \foreach \i/\x/\y in {
    0/0/0,1/2/0,2/4/0,3/6/0,4/8/0,5/10/0,6/12/0,
    7/2/2,8/2/4,9/8/2,10/6/-2,11/6/-4,12/10/-2}
    \coordinate (c\i) at (\x,\y);
  \foreach \i in {0,...,12}
    \filldraw[fill=black!12,draw=black!65]
      (c\i) circle (1);
  \foreach \i in {0,...,12}{
    \begin{scope}
      \foreach \k in {0,...,12}{
        \ifnum\k=\i\relax\else
          \pgfseteorule\clip
            (-3,-7) rectangle (15,7) (c\k) circle (2);
        \fi
      }
      \draw[dashed,black!40] (c\i) circle (2);
      \foreach \j in {0,...,12}{
        \ifnum\j=\i\relax\else
          \begin{scope}
            \clip (c\j) circle ({2+\gapwidth});
            \draw[thinred,line width=1.5pt] (c\i) circle (2);
          \end{scope}
        \fi
      }
    \end{scope}
  }
\end{tikzpicture}
\end{minipage}\hfill
\begin{minipage}[c]{.40\textwidth}
\centering
\begin{tikzpicture}[scale=.85]
  \definecolor{thinred}{RGB}{190,55,45}
  \definecolor{newgreen}{RGB}{40,124,104}
  \pgfmathsetmacro{\gapwidth}{0.18}
  \pgfmathsetmacro{\endangle}{acos((8-(2+\gapwidth)^2)/8)}
  \pgfmathsetmacro{\contactangle}{acos((8-(2+.65*\gapwidth)^2)/8)}
  \filldraw[fill=black!12,draw=black!65] (-1,0) circle (1);
  \filldraw[fill=black!12,draw=black!65] (1,0) circle (1);
  \draw[dashed,black!40] (-1,0) circle (2);
  \draw[dashed,black!40] (1,0) circle (2);
  \coordinate (incoming) at
    ({-1+2*cos(\contactangle)},{2*sin(\contactangle)});
  \draw[newgreen,line width=.8pt] (incoming) circle (1);
  \begin{scope}[thinred,line width=2pt,line cap=round]
    \draw ({-1+2*cos(60)},{2*sin(60)})
      arc (60:\endangle:2);
    \draw ({1+2*cos(180-\endangle)},{2*sin(180-\endangle)})
      arc ({180-\endangle}:120:2);
    \draw ({-1+2*cos(-\endangle)},{2*sin(-\endangle)})
      arc ({-\endangle}:-60:2);
    \draw ({1+2*cos(-120)},{2*sin(-120)})
      arc (-120:{-180+\endangle}:2);
  \end{scope}
  \fill[newgreen] (incoming) circle (1.5pt);
\end{tikzpicture}
\end{minipage}
\caption{Thin attachment regions for a branching chain of unit balls
(left) and a local configuration (right).}
\label{fig:thin-attachment}
\end{figure}

We next exclude compact shapes at optimal packing density.
Let $\omega_d\coloneqq|B(0,1)|$ 
and denote by $\Delta_d$ the optimal density of a packing of 
unit balls in $\R^d$.
For a bounded region $R$ of unit volume, define
\begin{equation}\label{dfn:PnR}
P_n(R)\coloneqq\left(\frac{\omega_d \, n}{\Delta_d}\right)^{1/d}R.
\end{equation} 
This normalization gives
$|P_n(R)|=\omega_d n/\Delta_d$ which is the volume that allows for
packing $n$ unit balls at optimal density. 

The following corollary shows that planar Euclidean \DLA{} leaves a
positive fraction of its centers outside every such region.
As in Corollary~\ref{cor:shapes}, the region
may be chosen after seeing the entire \DLA{} realization.

\begin{corollary}\label{cor:euclidean}
For $d=2$ there is a constant $\alpha>0$ such that almost surely, for every
compact set $R\subset\R^d$ of unit volume,
\[
 \liminf_{n\to\infty}\frac{\bigl|\{x\in C_n\mid x\notin P_n(R)\}\bigr|}{n}\ge\alpha.
\]
\end{corollary}

Rotational invariance of the \DLA{} evolution implies that any deterministic scaling limit---that is
also compact, connected and contains the origin---must be a ball.
Corollary~\ref{cor:euclidean} excludes the ball of asymptotically
smallest volume that could contain $n$ unit balls. The positive gap $\alpha$ also excludes balls
with slightly larger radii, while leaving open the possibility of a
sufficiently larger ball with a lower bulk density.
The deterministic packing comparisons needed to complete the argument 
(see Lemma \ref{eu-lem:localpacking}) are given here only for $d=2$, but 
should extend to $d \in \{3,8,24\}$ with more advanced results about optimal sphere packing
(see Remark \ref{rmk:sphere-packing}).

\subsection{Past results} \label{sec:past}
Kesten proved that the radius of \DLA{} on $\Z^d$ is almost
surely $O(n^{2/(d+1)})$ for $d=2$ and $d\ge4$,
and $O((n\log n)^{1/2})$ for $d=3$
\cite{Kesten1987Hitting,Kesten1987How,Kesten1990Upper}.
Lawler improved the three-dimensional bound to
$n^{1/2}(\log n)^{1/4}$
\cite[Section~2.6]{Lawler1991Intersections}.
Losev and Smirnov gave new proofs of the planar and
three-dimensional bounds and extended them to the
dielectric-breakdown model \cite{Losev2025How}.
Kesten also proved that, almost surely,
$|A_\infty\cap B(0,r)|\ge c_d r^{d-1}$ for arbitrarily
large $r$, where $c_d>0$ is deterministic
\cite{Kesten1990Upper}.
Benjamini and Yadin developed a general framework for
radius bounds that applies to a broad class of graphs
\cite{Benjamini2017Upper}.

For Euclidean \DLA{} in the plane, Bishop and Peres presented
an adaptation of Kesten's $n^{2/3}$ radius bound for aggregating
unit disks \cite[Theorem~1.6.1]{BishopPeresConformal}.
Losev and Smirnov also noted that their methods can be
adapted to unit particles in $\R^2$ and $\R^3$, although
their paper treats lattice models
\cite[Section~1]{Losev2025How}.

\subsection{Further questions}

We expect Theorem~\ref{thm:boundary} to remain true when
$\partial_e A_n$ is restricted to edges joining $A_n$ to
the infinite component of its complement. 
A related question concerns the external boundary of the
Euclidean exclusion region $\mathcal O_{C_n}$. Let $U_n$
be the unbounded component of
$\R^d\setminus\mathcal O_{C_n}$, and write
$\operatorname{SA}$ for surface area, interpreted as arc
length when $d=2$.
We conjecture that, for every $d\ge2$, almost surely
$\liminf_{n\to\infty}
 {\operatorname{SA}(\partial U_n)}/{n}>0.$
This would give linear growth of the external boundary
seen by the incoming particle's center.

The \DLA{} aggregate is believed to have a fractal scaling
limit \cite{Witten1981Diffusion}.
Consistent with this picture, we conjecture that, under
the assumptions of Corollary~\ref{cor:shapes}, almost surely
$|A_n\cap R_n|=o(n).$ Equivalently, the corollary should hold with $\alpha=1$.
One may also ask for an improvement over the
elementary lower bound on the radius by showing that 
almost surely,
$n^{-1/d} \max\{|x|\colon x\in A_n\} \to\infty$. 
The same question is natural for planar Euclidean \DLA{}.

\subsection{Organization}

Section~\ref{sec:proof-discussion} outlines the proofs of our results.
Starting with the lattice \DLA{} results, in Section~\ref{sec:boundary}
we prove Theorem~\ref{thm:boundary} using three lemmas, 
which are proved in Section~\ref{sec:lemmas}.
Section~\ref{sec:vacancies} proves the permanent vacancy estimate 
of Theorem~\ref{thm:vacancies}, while Section~\ref{sec:shapes} 
deduces Corollary~\ref{cor:shapes}.
For Euclidean \DLA{}, Section~\ref{sec:attachment} proves
Theorem~\ref{thm:attachment} with Section~\ref{sec:euclidean} then deriving
Corollary~\ref{cor:euclidean} 
by a deterministic packing comparison.

\section{Proof Ideas}\label{sec:proof-discussion}

\subsection{Lattice \DLA{}}

Our proof uses harmonic functions to relate the attachment rule
to the geometry of the aggregate. 
The functions are constructed
from a kernel $K$ given by a difference of discrete second
derivatives of the Green function or potential kernel. More
precisely, $$K=(2d)^{-1}(\partial_2^2-\partial_1^2)\Phi,$$
where $\partial_i^2$ denotes the 
discrete second derivative (see \eqref{def:kernel}).
Here $\Phi$ is the Green function in dimensions $d\ge3$ and the
negative potential kernel in dimension two. Taking the second differences 
defining $K$ improves the decay of the fundamental potential and gives the summability estimates
in Lemma~\ref{lem:kernel}.

A key feature of $K$ is that, for each $x\in A$, the translate $K(\,\cdot-x)$ is harmonic outside
$\overline A$. Thus, its value along the incoming
walk is a martingale until the walk first reaches a site adjacent
to $A$. 
Optional stopping implies the following identity (see Lemma~\ref{lem:cancellation}),
\begin{equation}\label{eq:mg-ident}
 \EE[K(X_{n+1}-x)\mid\cF_n]=0
 \qquad \forall x\in A_n,
\end{equation}
where $\cF_n\coloneqq\sigma(A_1,\ldots,A_n)$. 
We note that \eqref{eq:mg-ident} ensures that for any bounded $\mathcal F_k$-measurable coefficients $a_k(x)$, the process
\begin{equation}\label{eq:mg-K}
 \sum_{k=1}^{n-1}\sum_{x\in A_k}a_k(x)K(X_{k+1}-x) 
\end{equation}
is a martingale. This family may have other applications to \DLA{}.

To prove Theorem~\ref{thm:boundary}, we introduce an edge
energy function $F(A_n)$ which equals the sum of squares of a field
$\{\ell_e(A_n)\}$ over edges with both endpoints in $A_n$.
Specifically, for an edge $e=(x,z)$, define
\begin{equation}\label{dfn:field-n}
 \ell_e(A_n)\coloneqq
 \sum_{y\in A_n}\bigl(K(x-y)-K(z-y)\bigr).
\end{equation}
This measures how much the combined influence of the occupied
sites changes across the edge. 
We then use that convolution with
$K$ is a contraction on
$\ell^2(\Z^d)$ to deduce that $F(A) \leq |\partial_e A|$. 
Thus it suffices to prove that $F(A_n)$ grows linearly.

Subtracting the accumulated squared field increments from
$F(A_n)$ gives a submartingale $Q_n$ based on edge-indexed
processes analogous to those of \eqref{eq:mg-K}.
Lemmas~\ref{lem:kernel} and~\ref{lem:cancellation}
bound its increments, so Azuma--Hoeffding shows that
the energy is at least the accumulated squared increments
up to $o(n)$.
By Lemma~\ref{lem:energy}, attachments with an occupied
neighbor in direction $\pm e_1$ or $\pm e_2$ contribute
a uniformly positive squared increment.
If enough attachments have such a neighbor, the resulting
large energy forces a large boundary. Otherwise, many edges
in those directions must remain between $A_n$ and its
complement, again forcing a large boundary.
Corollary~\ref{cor:shapes} follows from Theorem~\ref{thm:boundary} since 
changing one lattice site changes the number of boundary edges by
at most $2d$. 

The proof of Theorem~\ref{thm:vacancies} closely follows that of
Theorem~\ref{thm:boundary}.  The main difference is that we localize
the energy to edges inside a fixed cube while retaining the field
\eqref{dfn:field-n} generated by the entire aggregate.
As before, attachments contribute uniformly positive squared
increments and spatial cancellation bounds the energy.
The cancellation ensures that the localized energy is small when the cube
is nearly full.
We also use the fact that each site is added at most once to bound
the sum of squared submartingale increments by a constant times
the cube's volume. A uniform over all
growth times concentration bound yields \eqref{eq:vacancy-tail}.

The cancellation identity \eqref{eq:mg-ident} is related to identities used in the random interlacements literature. Sznitman \cite[Equation~(1.8)]{Sznitman2010Vacant} records the constant-potential identity in dimensions $d\ge3$ and the planar counterpart is discussed by Comets, Popov, and Vachkovskaia \cite[Section~2.1]{Comets2016Two}. To our knowledge, the connection between this harmonic cancellation and an energy controlled by the boundary size and permanent vacancies is new.

\subsection{Euclidean \DLA{}}
An almost-double contact requires the incoming center to hit a small region 
on the sphere $\partial B(a,2)$ around the existing center $a \in C$ to which 
it attaches (see Figure~\ref{fig:thin-attachment}). In Lemma~\ref{eu-lem:thin}
we show that the surface area of this region is at most 
$O(\sqrt\eta)$ when $d=2$ and $O(\eta)$ when $d \ge 3$, 
uniformly in the $2$-separated set $C$ of existing centers, 
whereby the Poisson formula for the Brownian hitting density 
uniformly bounds the probability of such almost-double 
contact, starting anywhere on $\partial B(a,3)$. 
Theorem~\ref{thm:attachment}\textup{(a)}
then follows by showing that a constant multiple of
$\PP_z(\tau_C<\infty)$ bounds the expected number of such
radius-$3$ neighborhoods visited before $\tau_C$
(uniformly in $C$ and in starting points $z$ with $\dist(z,C)\ge3$).

To get Theorem \ref{thm:attachment}\textup{(b)}, note that each attachment creates exactly one near-contact edge, apart from having some universally bounded number 
of additional edges when almost-double contact occurs. From 
Theorem \ref{thm:attachment}\textup{(a)} the conditional expected number of new edges 
is arbitrarily close to one, so a standard martingale 
concentration bound implies that the total number of proximity edges 
must be close to the number of particles.

Corollary~\ref{cor:euclidean} follows because the low average degree
of the proximity graph in Theorem~\ref{thm:attachment}\textup{(b)}
is incompatible with nearly optimal packing in the plane.
Specifically, Lemma~\ref{eu-lem:localpacking} shows that average
degree at most five forces a fixed loss in packing density.
The proof uses a periodic tiling and the fact that a planar Voronoi
cell has minimum area only when its disk touches six others
\cite{Mazel2022Minimal}.

\section{Lattice \DLA{} Large Boundary: Proof of Theorem~\ref{thm:boundary}}\label{sec:boundary}
Before proving Theorem~\ref{thm:boundary} we state a few definitions and lemmas whose proofs we defer until
Section~\ref{sec:lemmas}.
We write $\ind_A$ for an indicator. Constants denoted by $c_d$ may change
from line to line. We begin by defining a kernel $K\colon\Z^d\to\R$ that is at the heart of our arguments. 

Let
\begin{equation}\label{dfn:semi-group}
 Pf(x)\coloneqq \frac{1}{2d}\sum_{i=1}^d
       \bigl(f(x+e_i)+f(x-e_i)\bigr)
\end{equation}
be the random walk averaging operator. A function $f$ is harmonic at
$x$ when $Pf(x)=f(x)$. Write $\delta_x\coloneqq\ind_{\{x\}}$ and
$p_k(x)\coloneqq\PP_0(S_k=x)$. Following  \cite[Sections~4.3--4.5]{Lawler2010Random}, for $d\ge3$, the random walk \emph{Green function} is
\[
 G(x)\coloneqq \sum_{k\ge0}p_k(x),
\]
and in dimension two the \emph{potential kernel} is
\[
 a(x)\coloneqq \lim_{N\to\infty}\sum_{k=0}^N\bigl(p_k(0)-p_k(x)\bigr).
\]
Define
\begin{equation}\label{def:kernel}
 K(x)\coloneqq \frac{1}{2d}\bigl(\Phi(x+e_2)+\Phi(x-e_2)
                         -\Phi(x+e_1)-\Phi(x-e_1)\bigr),
\end{equation}
where $\Phi\coloneqq G$ for $d\ge3$ and $\Phi\coloneqq -a$ for $d=2$. Writing $I$ for the identity
operator, we have $(I-P)\Phi=\delta_0$.
For $d\ge3$, $K(x)$ is $1/2d$ times the difference between the
expected numbers of visits to $\{x\pm e_2\}$ and $\{x\pm e_1\}$
by a simple random walk started at the origin. In dimension two, 
the same difference is obtained by counting
visits up to time $N$ and then letting $N\to\infty$. The factor $1/2d$ is included to make
convolution with $K$ a contraction.

Lemma~\ref{lem:kernel} records various kernel estimates that we shall need.
\begin{prooflemma}{proof:kernel}\label{lem:kernel}
The kernel $K$ is real-valued and even, $K(0)=0$, and
\begin{equation}\label{eq:source}
 (I-P)K=\frac{1}{2d}
  (\delta_{e_2}+\delta_{-e_2}-\delta_{e_1}-\delta_{-e_1}).
\end{equation}
This kernel has the tail bounds 
\begin{align}\label{eq:tails-K}
|K(x)|&\le c_d (1+|x|)^{-d}, \\
|K(x+e_i)-K(x)|&\le c_d(1+|x|)^{-d-1} \,, \qquad 1 \le i \le d \,.
\label{eq:K-diff-tails}
\end{align}
Moreover, the linear operator
\begin{equation}\label{def:T-operator}
 Tf(x)\coloneqq \sum_y K(x-y)f(y),
\end{equation}
initially defined for finitely supported $f$, extends uniquely
to a contraction on $\ell^2(\Z^d)$
\begin{equation}\label{eq:T-contraction}
 \|Tf\|_2\le\|f\|_2
 \qquad\text{for every }f\in\ell^2(\Z^d).
\end{equation} Lastly, 
\begin{equation}\label{def:kappa}
 \kappa_d\coloneqq |K(e_1)|=|K(e_2)|
 >\kappa_d^*\coloneqq\max_{|z|_1=2}|K(z)|.
\end{equation}
\end{prooflemma}

For concreteness, orient each nearest-neighbor edge from the endpoint closer
to the origin to the endpoint farther from the origin.
Any other choice of orientations would work equally well.
For an edge $e$ oriented from $x$ to $z$, any $y\in\Z^d$,
and any finite set $A$, define
\begin{equation}\label{dfn:field}
 \ell_e(y)\coloneqq K(y-x)-K(y-z),
 \;\;
 \ell_e(A)\coloneqq\sum_{y\in A}\ell_e(y)
 =(T\ind_A)(x)-(T\ind_A)(z).
\end{equation}
Since $K$ is a difference of second discrete derivatives of $\Phi$,
$\ell_e(y)$ takes one further difference across $e$, while
$\ell_e(A)$ sums these field differences over $y\in A$.

Let $\mathcal{E}(A)$ denote the set of edges with
both endpoints in $A$. Define the \emph{edge energy} $F(A)$
and the \emph{increment energy} $W_A(y)$ by
\begin{equation}\label{dfn:bdry-energy}
 F(A)\coloneqq\sum_{e\in \mathcal{E}(A)}\ell_e(A)^2,
 \qquad
 W_A(y)\coloneqq\sum_{e\in \mathcal{E}(A)}\ell_e(y)^2.
\end{equation}

\begin{prooflemma}{proof:energy}\label{lem:energy}
Every finite set $A\subset\Z^d$ satisfies
\begin{equation}\label{eq:bdry-energy-ineq}
 |\partial_e A|\ge F(A).
\end{equation}
If $A$ is connected, $|A|\ge 2$, and $y\notin A$ has a neighbor
in $A$ along either the $e_1$-axis or the $e_2$-axis, then
\begin{equation}\label{eq:W-increment}
 W_A(y)\ge v_d\coloneqq(\kappa_d-\kappa_d^*)^2>0.
\end{equation}
\end{prooflemma}

The next lemma uses that translates of $K$ are harmonic away from the aggregate to obtain a martingale related to the next attachment. 
We then derive the consequence of this for changes in $F$ and $W$. To state the result, in view of \eqref{eq:K-diff-tails}, we define the constant
\begin{equation}\label{eq:edge-summ}
J_d \coloneqq\sum_{i=1}^d\sum_{x\in\Z^d} |K(x+e_i)-K(x)|<\infty.
\end{equation}
\begin{prooflemma}{proof:cancellation}\label{lem:cancellation}
Let $A$ be a finite nonempty connected set and $X$ a site on $\partial A$ 
chosen according to the harmonic measure $\mu_A$. For every $x\in A$,
\begin{equation}\label{eq:K-cancellation}
 \EE[K(X-x)\mid A]=0.
\end{equation}
Moreover, 
\begin{equation}\label{dfn:DAX}
 D(A,X)\coloneqq F(A\cup\{X\})-F(A)-W_A(X), 
\end{equation}
satisfies 
\begin{equation}\label{eq:sub-MG}
 \EE[D(A,X)\mid A]\ge0\quad {\rm and} \quad
|D(A,X)|\le (2+2d)J_d^2 =: L_d \,.
\end{equation}
\end{prooflemma}

\begin{proof}[Proof of Theorem~\ref{thm:boundary}] Define the $\cF_n$ adapted sequences
\begin{equation}\label{def:squared-increments}
 V_n\coloneqq\sum_{k=2}^n W_{A_{k-1}}(X_k),
 \qquad Q_n\coloneqq F(A_n)-V_n
\end{equation}
($V_n$ records the squared increments only at edges already
occupied when the new particle arrives). We first prove that
\begin{equation}\label{eq:edge-lower-tail}
 \PP(Q_n\le-\varepsilon n)
 \le\exp\{-\varepsilon^2 n/(8L_d^2)\},
 \qquad \varepsilon>0.
\end{equation} 
To this end, note that $Q_1=0$ and for any $k \ge 2$, 
\begin{equation}\label{eq:subMG-diff}
D_k \coloneqq Q_k-Q_{k-1} = F(A_k)-F(A_{k-1})-W_{A_{k-1}}(X_k) = D(A_{k-1},X_k),
\end{equation}
where the last equality uses $A_k=A_{k-1}\cup\{X_k\}$ and \eqref{dfn:DAX}.
By construction,
\[
M_n \coloneqq Q_n - \sum_{k=2}^n \EE[ D_k \mid\cF_{k-1}] \le Q_n
\]
(see the \abbr{lhs} of \eqref{eq:sub-MG}), 
is an $(\cF_n)$ martingale, starting at $M_1=Q_1=0$, 
whose increments are bounded by $2L_d$
(see the \abbr{rhs} of \eqref{eq:sub-MG}). Applying the Azuma--Hoeffding
inequality to this martingale
gives \eqref{eq:edge-lower-tail}.

Now, let $H(A)$ count the number of edges in $\mathcal{E}(A)$ which are parallel to $e_1$ or $e_2$. Since $H(A_{k+1})-H(A_k) \le 4$ for any $k$, at least $(H(A_n)-2)/4$ 
of our attachments for $3 \le k \le n$ must add to the \DLA{} aggregate 
some edge which is parallel to $e_1$ or $e_2$. By \eqref{eq:W-increment} of Lemma~\ref{lem:energy}, each such attachment contributes at least $v_d$ to $V_n$, hence
\begin{equation}\label{eq:lbd-Vn}
V_n \ge \frac{v_d}{4} (H(A_n)-2) \,.
\end{equation}
Further, from \eqref{eq:edge-lower-tail} we deduce by Borel--Cantelli that 
\begin{equation}\label{eq:energy-growth}
 F(A_n)\ge V_n-o(n) \qquad\text{almost surely}.
\end{equation}
If $H(A_n) \ge n$ then by these two facts and the assertion 
\eqref{eq:bdry-energy-ineq} of Lemma~\ref{lem:energy}, 
\[
 |\partial_e A_n|\ge F(A_n)\ge V_n-o(n)
 \ge \frac{v_d}{4}n-o(n).
\]
Otherwise, when $H(A_n)<n$, count the four edges incident to each site
of $A_n$ in the directions $\pm e_1$ and $\pm e_2$.
Each internal edge contributes twice to this count,
and each boundary edge contributes once, so 
\[
 |\partial_e A_n|\ge4n-2H(A_n)>2n.
\]
Taking $b_d\coloneqq 2 \wedge \frac{v_d}{4} >0$, 
our two lower bounds on $|\partial_e A_n|$ complete the proof.
\end{proof}

\begin{remark}
For every $0<c<b_d$, the deterministic bounds above and
\eqref{eq:edge-lower-tail} give $b=b(d,c)>0$ such that
\[
 \PP\bigl(|\partial_e A_n|\le cn\bigr)\le e^{-bn}
\]
for all sufficiently large $n$.
\end{remark}

\section{Proofs of the lattice lemmas}\label{sec:lemmas}

We proceed to prove the analytic estimates of
Lemma~\ref{lem:kernel}, the deterministic energy bounds of 
Lemma~\ref{lem:energy} and the cancellation and bounded-increment
properties of Lemma~\ref{lem:cancellation}, in that order. 

\begin{proof}[Proof of Lemma~\ref{lem:kernel}]\hypertarget{proof:kernel}{}
Starting with \eqref{eq:source}, recall that $(I-P)\Phi=\delta_0$
for the random walk averaging operator $P$ of \eqref{dfn:semi-group},
our choice of $\Phi$ in \eqref{def:kernel} and any $d \ge 2$. The identity 
\eqref{eq:source} follows, for $(I-P)$ is linear and commutes with translations.
The kernel $K$ is real-valued and even, since $\Phi$ is real-valued
and even. Moreover, exchanging coordinates $1$ and $2$ reverses the sign of
$K$, so $K(0)=0$.

For decay, the standard asymptotics
\cite[Theorems~4.3.1 and~4.4.4]{Lawler2010Random} imply that
\begin{align}\label{eq:L-asympt}
 G(x)&=a_d|x|^{2-d}+O_d(|x|^{-d}), \qquad 
 a(x)=\frac{2}{\pi}\log|x|+C_0+O(|x|^{-2}).
\end{align}
Write $g(x)=a_d|x|^{2-d}$ for $d\ge3$ and
$g(x)=-(2/\pi)\log|x|-C_0$ for $d=2$.
The bounds at \eqref{eq:L-asympt} imply that $\Phi(x)=g(x)+O_d(|x|^{-d})$.
The second discrete differences of the smooth function $g$
are $O_d(|x|^{-d})$. This follows by integrating its second
derivatives along unit segments away from the origin.
Substituting into \eqref{def:kernel}, and taking $c_d$ large enough to handle bounded $x$, therefore gives
\eqref{eq:tails-K}.

Turning to \eqref{eq:K-diff-tails}, let
$r=\lfloor |x|/4\rfloor$ for sufficiently large $|x|$.
By \eqref{eq:source}, the kernel $K$ is harmonic in the lattice ball
of radius $r$ centered at $x$.
The discrete difference estimate of 
\cite[Theorem~6.3.8(a)]{Lawler2010Random}, together with
\eqref{eq:tails-K}, gives that
\[
 |K(x+e_i)-K(x)|
 \le \frac{c_d}{r}
    \sup_{\substack{z\in\Z^d\\|z-x|\le r+1}}|K(z)|
 \le c_d|x|^{-d-1}
\]
and increasing $c_d$ handles \eqref{eq:K-diff-tails} for 
the remaining, bounded, values of $x$.

Next, for finitely supported $f$, write
$\widehat f(\theta)\coloneqq \sum_xf(x)e^{-ix\cdot\theta}$.
Plancherel extends this transform to an isometry from
$\ell^2(\Z^d)$ onto $L^2([-\pi,\pi]^d,(2\pi)^{-d}d\theta)$.
Now, from \eqref{eq:tails-K} we know that $K\in\ell^2(\Z^d)$ and 
setting the non-negative $\alpha_i(\theta)\coloneqq 1-\cos\theta_i$, 
the linear operator $P$ of \eqref{dfn:semi-group} is 
bounded on $\ell^2(\Z^d)$, corresponding to 
the multiplier $\widehat P(\theta)\coloneqq 1 - d^{-1}\sum_{k} \alpha_k(\theta)$ 
on $L^2$. Thus, taking the Fourier transform of \eqref{eq:source}, now 
as an identity in $\ell^2(\Z^d)$, gives
\[
 (1-\widehat P(\theta))\widehat K(\theta)
 =d^{-1} (\alpha_1(\theta)-\alpha_2(\theta)),
\]
from which it follows that $\widehat K=m$ almost everywhere, where
\begin{equation}\label{def:K-Fourier}
 m(\theta)\coloneqq 
 \frac{\alpha_1(\theta)-\alpha_2(\theta)}{\sum_{k=1}^d\alpha_k(\theta)},
 \qquad 
 \quad(\theta\ne0),\qquad m(0)\coloneqq 0.
\end{equation}
While $m(\theta)$ is discontinuous at
zero (with limits $1$ and $-1$ along the first and second coordinate
axes), its value at that one point has no effect on the corresponding $L^2$ operator.
Acting on finitely supported functions $f$, the convolution $Tf$ in \eqref{def:T-operator}
has the Fourier transform $m\widehat f$. The bound $|m|\le1$ which is immediate from 
\eqref{def:K-Fourier}, implies by Plancherel that
$\|Tf\|_2\le\|f\|_2$ for all such functions. With the
linear subspace of finitely supported functions dense in $\ell^2(\Z^d)$, we thus deduce \eqref{eq:T-contraction}
that $T$ extends uniquely to a contraction on $\ell^2(\Z^d)$.

Write $\langle\cdot\rangle$ for normalized integration over the
Fourier torus. The interchange of coordinates $1$ and $2$ reverses the
signs of $K$ and $m$. Thus, $\langle m\rangle=0$ and by \eqref{def:K-Fourier}
\[
-K(e_1) = K(e_2) 
 =-\left\langle
       \frac{\alpha_2(\alpha_1-\alpha_2)}{\sum_k\alpha_k}
   \right\rangle
 =\frac{1}{2}\left\langle
       \frac{(\alpha_1-\alpha_2)^2}{\sum_k\alpha_k}
   \right\rangle = \kappa_d > 0 .
\]
Equation~\eqref{eq:source} shows that $K$ is harmonic away
from $\{\pm e_1,\pm e_2\}$, while \eqref{eq:tails-K} shows
that it vanishes at infinity. Its absolute value at each of these four source
sites is $\kappa_d$. The maximum principle therefore gives
$|K|\le\kappa_d$ everywhere. Every site at lattice distance two
from the origin belongs to the unbounded component of the complement
of the four source sites. The strong maximum principle, together with
decay at infinity, gives $|K(z)|<\kappa_d$ at each such site. 
Define $\kappa_d^*\coloneqq\max_{|z|_1=2}|K(z)|$. The maximum is
over a finite set, so $\kappa_d^*<\kappa_d$.
\end{proof}

\begin{proof}[Proof of Lemma~\ref{lem:energy}]\hypertarget{proof:energy}{}
The convolution $T$ with an even kernel commutes with taking lattice differences
and is a contraction on $\ell^2(\Z^d)$ by \eqref{eq:T-contraction}. Hence, 
\begin{align*}
\sum_{i=1}^d \sum_{x \in \Z^d}
       |(T\ind_A)(x+e_i)-(T\ind_A) (x)|^2
 &=\sum_{i=1}^d
       \|T(\ind_A(\,\cdot+e_i)-\ind_A(\cdot))\|_2^2\\
 &\le \sum_{i=1}^d
       \|\ind_A(\,\cdot+e_i)-\ind_A (\cdot) \|_2^2
 =|\partial_e A|.
\end{align*}
Further, by \eqref{dfn:bdry-energy} we have that
\begin{align}\label{eq:ubd-energy}
 F(A) \le \sum_{e \in \mathcal{E}(\Z^d)} \ell_e(A)^2
 & = \sum_{i=1}^d \sum_{x \in \Z^d}
       |(T\ind_A)(x+e_i)-(T\ind_A) (x)|^2 \le |\partial_e A|,
\end{align} 
which thereby proves the assertion \eqref{eq:bdry-energy-ineq} for every finite $A$.

Turning to show \eqref{eq:W-increment}, choose $x\in A$ with
$y-x\in\{\pm e_1,\pm e_2\}$. With $A$ connected and $|A| \ge 2$, there must be some 
$\{x,z\} \in \mathcal{E}(A)$. Moreover, since $y \not\in A$ necessarily $z\ne y$ with 
$|y-z|_1=2$. Consequently, thanks to \eqref{def:kappa}, 
\[
 |\ell_{\{x,z\}}(y)|
 \ge |K(y-x)|-|K(y-z)|
 \ge\kappa_d-\kappa_d^*.
\]
This edge contributes at least $v_d$ to $W_A(y)$, and since all other edge contributions 
are non-negative, our second assertion \eqref{eq:W-increment} follows.
\end{proof}

\begin{proof}[Proof of Lemma~\ref{lem:cancellation}]\hypertarget{proof:cancellation}{}
Starting with the proof of \eqref{eq:K-cancellation}, fix $x\in A$ and start an unconditioned walk at
$z\notin\overline A$. The four sites $x\pm e_1,x\pm e_2$ belong to
$\overline A$. At every time $k<\tau_A$, the walk is outside
$\overline A$, so \eqref{eq:source} gives
$PK(S_k-x)=K(S_k-x)$. Thus,
$K(S_{k\wedge\tau_A}-x)$ is a bounded martingale, and for every
deterministic $k$,
\begin{equation}\label{eq:MG-stopped}
 \EE_zK(S_{k\wedge\tau_A}-x)=K(z-x).
\end{equation}
In case $d=2$, by the recurrence of the walk, $\tau_A<\infty$ almost surely.
For $d\ge3$, transience gives $|S_k|\to\infty$ almost surely, so 
$|K(S_{k \wedge \tau_A}-x)| \to 0$ on $\{\tau_A=\infty\}$. By bounded convergence, 
it thus follows from \eqref{eq:MG-stopped} that for any $d \ge 2$, $x \in A$ and $z \notin \overline A$
\begin{equation}\label{eq:stopped-K}
 \EE_z[K(S_{\tau_A}-x)\ind_{\{\tau_A<\infty\}}]=K(z-x).
\end{equation}
For $d=2$, when taking $|z|\to\infty$ we have by \eqref{eq:tails-K} 
that $|K(z-x)| \to 0$ and the convergence to $\mu_A$ of the law of $S_{\tau_A}$ 
(see \eqref{def:lattice-harmonic-measure}) yields \eqref{eq:K-cancellation}. 
For $d\ge3$, note that the first hitting time $\tau_x$ of $x \in \overline A$, 
must exceed $\tau_A$. Hence, for sufficiently large $|z|$, 
\[
 \PP_z(\tau_A<\infty)\ge\PP_z(\tau_x<\infty)
 =\frac{G(z-x)}{G(0)}\ge c_{d,x}|z|^{2-d}
\]
using the Green function asymptotic in
\cite[Theorem~4.3.1]{Lawler2010Random}. It thus follows upon dividing 
\eqref{eq:stopped-K} by $\PP_z(\tau_A<\infty)$ and using \eqref{eq:tails-K}, that as $|z| \to \infty$,
\[
 \left|\EE_z[K(S_{\tau_A}-x)\mid\tau_A<\infty]\right|
 \le C_{d,x}|z|^{-2}\longrightarrow 0.
\]
To get \eqref{eq:K-cancellation}, note that 
the expectation on the \abbr{lhs} converges to $\EE[K(X-x)\mid A]$
(by the convergence \eqref{def:lattice-harmonic-measure}
of the conditional law of $S_{\tau_A}$ to that of $X$, both supported on 
the fixed finite set $\partial A$). 

Turning to prove \eqref{eq:sub-MG}, for any $X \not\in A$,
\[
\ell_e (A \cup \{X\}) = \ell_e (A) + \ell_e (X),
\]
so expanding the terms $\ell_e(A \cup \{X\})^2$ in \eqref{dfn:bdry-energy}
gives that 
\begin{align}
 D(A,X) =F(A \cup \{X\}) - F(A) - W_A(X)  
 &=2\sum_{e\in \mathcal{E}(A)}
                   \ell_e (A) \,\ell_e(X)\notag\\
 &\quad+
       \sum_{e\in \mathcal{E}(A \cup \{X\} )\setminus \mathcal{E}(A)}
                   \ell_e (A \cup\{X\})^2.
 \label{eq:edge-increment}
\end{align}
For an old edge $e=\{x,z\} \in \mathcal{E}(A)$, considering \eqref{eq:K-cancellation} 
at both endpoints $x,z$ gives 
\begin{equation}\label{eq:field-cancellation}
 \EE[\ell_e(X)\mid A]
 =\EE[K(X-x)-K(X-z)\mid A]=0,
\end{equation}
whereas the second term on the \abbr{rhs} of 
\eqref{eq:edge-increment} is nonnegative, so
$\EE[D(A,X)\mid A]\ge0$.

Further, orienting, without loss of generality, each lattice edge in a positive
coordinate direction, we obtain from \eqref{eq:edge-summ} that for 
\begin{equation}\label{dfn:hB}
h_B(y) \coloneqq \sum_{e\in \mathcal{E}(B)}|\ell_e(y)|
\end{equation}
and any finite set $B$,
\begin{align}
\max_{y \in \Z^d} h_B(y) &\le 
\max_{y \in \Z^d} \sum_{i=1}^d \sum_{z \in \Z^d} |K(y-z+e_i)-K(y-z)| = J_d,
\label{eq:unif-hB}
\\
\label{eq:bd-ell-A}
\max_{e \in \mathcal{E}(\Z^d)} |\ell_e(B)| &\le \max_{z \in \Z^d} 
\sum_{i=1}^d \sum_{y \in \Z^d} |K(y-z+e_i)-K(y-z)| = J_d.
\end{align}
At most $2d$ new edges appear at one attachment, so from \eqref{eq:bd-ell-A} 
we deduce that the second term in \eqref{eq:edge-increment} is at most $2dJ_d^2$, 
while the first term in \eqref{eq:edge-increment} has 
absolute value at most $2 J_d h_A(X)$. Taking into consideration 
also \eqref{eq:unif-hB} gives 
\begin{equation}\label{eq:bd-DAX}
|D(A,X)|\le 2 J_d h_A(X) + 2 d J_d^2 \le L_d,
\end{equation}
as stated in \eqref{eq:sub-MG}.
\end{proof}

\section{Lattice \DLA{} Permanent Vacancies: Proof of Theorem~\ref{thm:vacancies}}\label{sec:vacancies} 

We follow the proof of Theorem~\ref{thm:boundary}, utilizing
the field $\ell_e(A)$ of \eqref{dfn:field}, but with all edge sums restricted to
$\mathcal E(B)$.
Namely, we set here
\begin{align}\label{dfn:bdry-energy-B}
 F_B (A)&\coloneqq\sum_{e\in\mathcal E(A \cap B)} \ell_e(A)^2,
 \qquad
 W_{A \cap B} (y)\coloneqq\sum_{e\in\mathcal E(A \cap B)}\ell_e(y)^2,\\
\label{dfn:DAX-B}
 D_B (A,X)&\coloneqq F_B(A\cup\{X\})-F_B(A)-W_{A \cap B} (X).
\end{align}

Restricting the expansion in \eqref{eq:edge-increment} to
$\mathcal E(B)$ and utilizing the cancellation in \eqref{eq:field-cancellation}
recovers the \abbr{lhs} of \eqref{eq:sub-MG} for $D_B(A,X)$. 
Hence,
the absolute value of the first term on the \abbr{rhs} of that 
expansion is at most $2 J_d h_B(X)$ for $h_B(\cdot)$ of \eqref{dfn:hB},
and the second term in that expansion is zero whenever $X\notin B$.
We can therefore replace \eqref{eq:bd-DAX} by the refined bound
\begin{equation}\label{eq:bd-DAX-b}
|D_B(A,X)| \le g_B(X) \le L_d \quad {\rm where} \quad
g_B(y)\coloneqq 2J_d h_B(y)+2dJ_d^2\ind_B(y) \,.
\end{equation}
Moreover,
it follows from
\eqref{eq:edge-summ} that
\begin{equation}\label{eq:bd-sum-hB}
\sum_{y \in \Z^d} h_B(y)
=\sum_{e\in\mathcal E(B)}\sum_{y\in\Z^d} |\ell_e(y)|
\le {|\mathcal E(B)|} J_d \le d \, J_d |B|.
\end{equation}

Armed with the bounds \eqref{eq:bd-DAX-b} and \eqref{eq:bd-sum-hB},
we set the corresponding submartingale 
\begin{equation}\label{def:squared-increments-B}
Q_n\coloneqq F_B(A_n)-V_n \quad {\rm where} \quad 
 V_n \coloneqq\sum_{k=2}^n W_{A_{k-1} \cap B} (X_k)
\end{equation}
and fixing $a \in (0,1)$,  
proceed to show that for some $\beta=\beta(d) >0$,
\begin{equation}\label{eq:vacancy-local-tail}
 \PP\bigl( \inf_{k \ge 1} \{ Q_k \} \le -a |B| \bigr)
 \le e^{-\beta a^2 |B|}.
\end{equation}
Indeed, 
as in the proof of Theorem~\ref{thm:boundary} (see \eqref{eq:subMG-diff}), 
the differences 
\[
D_k \coloneqq Q_k-Q_{k-1} 
= D_B(A_{k-1},X_k)
\]
are such that $\EE[D_k\mid\cF_{k-1}]\ge 0$. Moreover, 
in view of \eqref{eq:bd-DAX-b}-\eqref{eq:bd-sum-hB},
for some $c_d$ finite and any finite set $B$, 
\begin{equation}\label{eq:bd-quad-var}
|D_k| \le L_d \quad {\rm and} \quad
\sum_{k \ge 2} D_k^2 \le L_d \sum_{k \ge 2} g_B(X_k) \le c_d |B|
\end{equation}
(since each site is attached to the aggregate at most once).

Unlike in the proof of Theorem~\ref{thm:boundary}, we need a concentration bound uniform over all time. 
For $0<\lambda\le(2L_d)^{-1}$, the inequality
$e^{-t-t^2}\le1-t$ for $|t|\le1/2$, shows that the non-negative 
\[
Z_n \coloneqq \exp\big(-\lambda Q_n-\lambda^2\sum_{k=2}^n D_k^2\big) = 
e^{ -\lambda D_n - \lambda^2 D_n^2} \, Z_{n-1}
\]
is a supermartingale starting at $Z_1=1$.
Choosing $\lambda=\min\{a/(2c_d),1/(2L_d)\}$ and $\beta=1/(4c_d \vee 4L_d)$
guarantees that $\lambda a - \lambda^2 c_d \ge \lambda a/2 \ge \beta a^2$.
Thanks to \eqref{eq:bd-quad-var} we then get by applying Doob's maximal inequality 
to $(Z_n)$ that 
\begin{align}\label{eq:Q-max}
\PP\bigl(\inf_{k \ge 1} \{  Q_k \} \le -a |B| \bigr) \le 
\PP\bigl( \sup_{k \ge 1} \{  Z_k \} \ge e^{\lambda a |B| - \lambda^2 c_d |B|}\big) 
\le e^{-\beta a^2|B|} \,,
\end{align}
as stated in \eqref{eq:vacancy-local-tail}.

We next bound \(F_B(A_n)\) from above and \(V_n\) from below to show that $A_\infty$ nearly filling \(B\) forces \(Q_n\) below a negative multiple of \(|B|\) thus allowing us to apply \eqref{eq:vacancy-local-tail}.
To this end, since $K(y)=K(-y)$ has
absolutely summable first differences 
(see \eqref{eq:K-diff-tails} or the \abbr{rhs} of \eqref{eq:bd-ell-A}),
we have for $e=\{x,z\}$ and any finite set $A$,
\[
\ell_e(A)+\ell_e(\Z^d \setminus A) = \sum_{y \in \Z^d} \ell_e(y) = 0
\]
(pairing $y$ with its reflection $(x+z)-y$ 
through the midpoint of $e$). Consequently, 
\begin{equation}\label{eq:ubd-energy-B}
F_B(A) = \sum_{e\in\mathcal E(A \cap B)} \ell_e(\Z^d \setminus A)^2 
\le  2 \sum_{e \in \mathcal E(B)} \ell_e(B^c \setminus A)^2 + 
2 \sum_{e \in \mathcal E (\Z^d)} \ell_e(B \setminus A)^2 \,.
\end{equation}

To control the first sum on the \abbr{rhs}, note that, by
\eqref{dfn:field} and \eqref{eq:K-diff-tails}, for some finite $c_d$
and all $e=\{x,z\}\in\mathcal E(B)$,
\[
\ell_e(B^c \setminus A)^2
 \le \Big( \sum_{y \in B^c} |K(y-x)-K(y-z)| \Big)^2 \le 
\frac{c_d}{4} [ \dist_\infty(x,B^c)^{-2} 
+ \dist_\infty(z,B^c)^{-2} ] \,.
\]
Consequently, 
\begin{equation}\label{dfm:s-B}
\sum_{e \in \mathcal E(B)} \ell_e(B^c \setminus A)^2 \le c_d \, s(B),  \qquad 
{\rm where} \qquad 
 s(B) \coloneqq\sum_{v\in B} \dist_\infty(v,B^c)^{-2}  
 \,.
\end{equation}
As in \eqref{eq:ubd-energy}, we bound the right-most  
sum in \eqref{eq:ubd-energy-B}
by $|\partial_e (B \setminus A)|$, to arrive at  
\begin{equation}\label{eq:vacancy-energy-bound}
 F_B(A) \le 2 |\partial_e (B \setminus A)|+ c_d \, s(B). 
\end{equation}

We next derive the following localized version of \eqref{eq:lbd-Vn},
\begin{equation}\label{eq:lbd-Vn-loc}
 V_n\ge\frac{v_d}{4} \bigl(H(A_n\cap B)- 4 s(B) -1\bigr),
\end{equation}
where $H(A\cap B)$ counts the edges in $\mathcal E(A\cap B)$ 
parallel to $e_1$ or $e_2$. Indeed, at most 
$4|\partial_{\mathrm{in}}B| \le 4 s(B)$ of such edges 
have an endpoint in 
\[
 \partial_{\mathrm{in}} B
  \coloneqq \{v \in B \colon \dist_\infty(v,B^c) = 1 \} 
\]
and when any other edge $\{x,X_k\}$ counted by $H(A_n\cap B)$ with  $3 \le k \le n$ is 
created, by the connectivity of $A_{k-1}$ the point
$x \in (A_{k-1} \cap B) \setminus \partial_{\mathrm{in}}B$ 
has a neighbor $z \in A_{k-1} \cap B$. While 
$A_{k-1} \cap B$ may not be connected, applying the proof of \eqref{eq:W-increment}
to the edge $\{x,z\} \in \mathcal{E}(A_{k-1} \cap B)$, gives 
$W_{A_{k-1}\cap B}(X_k)\ge v_d$. Thus, excluding the at most one such edge created at $k=2$
and noting that at most four edges counted by $H(A_n\cap B)$
are created at each attachment gives \eqref{eq:lbd-Vn-loc}.

Combining \eqref{eq:vacancy-energy-bound} and \eqref{eq:lbd-Vn-loc} results in 
\begin{equation}\label{eq:ubd-QnB}
 Q_n=F_B(A_n)-V_n
 \le 2 |\partial_e (B \setminus A_n)| 
 -\frac{v_d}{4} (H(A_n \cap B)-1) + (c_d + v_d) \, s(B) \,.
\end{equation}
(a). Turning to show \eqref{eq:vacancy-tail}, set
$N_n\coloneqq|B\setminus A_n|$ and note that 
$|\partial_e (B \setminus A_n)| \le 2 d N_n$. Moreover, 
$H(A_n\cap B)\ge H(B)-4N_n$ (since any site missing 
from $B$ removes at most four edges from $H(B)$). Plugging this in 
\eqref{eq:ubd-QnB} yields
\begin{equation}\label{eq:ubd-QnB-alpha}
 Q_n
 \le(4d+v_d) N_n-\frac{v_d}{4} (H(B)-1) + (c_d + v_d) \, s(B) \,.
\end{equation}
For some $\alpha>0$, any lattice cube of side $m \ge 1/\alpha$ 
has $H(B) \ge |B| + 2$ and  
\begin{equation}\label{eq:sB-ubd}
s(B) \le 
2 \, d \, m^{d-1}\sum_{i\ge 1} i^{-2} < \frac{5 d}{m}  |B|\le
5 d \, \alpha |B| \,.  
\end{equation}
Reducing $\alpha$ such that also 
\[
(4d+v_d+ 5 d(c_d + v_d)) \alpha \le \frac{v_d}{8} \quad {\rm and } \quad 
\alpha \le \beta (v_d/8)^2 \,,
\]
we then deduce from \eqref{eq:ubd-QnB-alpha} that $N_n\le\alpha|B|$ 
implies $Q_n\le-v_d|B|/8$. Hence, by \eqref{eq:vacancy-local-tail}
\[
 \PP (N_n \le \alpha |B| ) \le  
 \PP( Q_n \le - \frac{v_d}{8} |B| ) \le e^{- \alpha |B|} \,,
\]
which since $N_n \downarrow |B \setminus A_\infty|$, concludes the proof of 
\eqref{eq:vacancy-tail}.

\noindent
(b). To prove \eqref{dfn:loc-edge-bdry}, we set $U_n \coloneqq A_n \cap B$ and 
note that 
\[
|\partial_e(B \setminus A_n)| \le |\partial_e U_n| + |\partial_e B|
\le |\partial_e U_n| + 2 d |\partial_{\mathrm{in}} B| \le 
|\partial_e U_n| + 2 d s(B) \,,
\]
while counting edges parallel to $e_1$ or $e_2$ that touch $U_n$ results in 
\[
4|U_n|\le 2H(U_n)+|\partial_e U_n|.
\]
Plugging those bounds in \eqref{eq:ubd-QnB} gives 
for some $c'_d$ finite 
and $\gamma=2v_d/(16+v_d)$, that if $|\partial_e U_n| \le \gamma |U_n|$ then 
\[
Q_n
 \le (2 + \frac{v_d}{8} ) |\partial_e U_n| 
 -\frac{v_d}{2} |U_n| + \frac{c_d'}{5d} \, s(B) 
\le - \frac{v_d}{4} |U_n| + \frac{c_d'}{5d} \, s(B) \,. 
\]
Recall \eqref{eq:sB-ubd}
that $s(B) \le 5 d \, \gamma \rho |B|$ for any $B$ of side at least $1/(\gamma \rho)$.
Further reducing $\gamma \le v_d/(8 c_d')$ we consequently get by the preceding, 
that if also $|U_n| \ge \rho |B|$ then 
\[
Q_n \le  - \Big(\frac{v_d}{4} - c_d' \gamma\Big) \rho |B| \le - \frac{v_d}{8} \rho |B| \,.
\]
Taking also $\gamma \le \beta (v_d/8)^2$ we thus deduce from 
\eqref{eq:vacancy-local-tail} that
\[
\PP\Bigl(
 |B\cap A_n|\ge \rho|B|,
 \quad
 |\partial_e(B\cap A_n)|\le\gamma |B \cap A_n|
 \Bigr) \le \PP(Q_n \le - v_d \rho |B|/8) \le e^{-\gamma \rho^2 |B|}.
\]
Taking $n \to \infty$ yields \eqref{dfn:loc-edge-bdry} 
since the $\Z_+^2$-valued variables 
$(|B\cap A_n|,|\partial_e(B\cap A_n)|)$ converge in law to 
$(|B\cap A_\infty|,|\partial_e(B\cap A_\infty)|)$.

\begin{remark}\label{rmk:general-boxes} 
The preceding proof, after adjusting constants, extends to collections 
of finite sets $B \subset \Z^d$ indexed by $\alpha>0$, for which
$\inf_B \{ H(B)/|B| \}$ is uniformly bounded away from zero while  
$\sup_B \{ s(B)/|B| \} \le \alpha$.
\end{remark}

\section{Lattice \DLA{} Shape Exclusion: Proof of Corollary~\ref{cor:shapes}}\label{sec:shapes}

For any finite lattice sets $A,B$,
\begin{equation}\label{eq:boundary-stability}
 \bigl||\partial_e A|-|\partial_e B|\bigr|
 \le 2d|A\triangle B|.
\end{equation}
Indeed, an edge belonging to exactly one of the two boundaries has
an endpoint in $A\triangle B$, and each vertex meets at most $2d$
edges. Since $|A_n|=n$, we also have
\[
 |A_n\triangle B|=2|A_n\setminus B|+|B|-n.
\]
Consequently,
\begin{equation}\label{eq:shape-boundary-comparison}
 \frac{|A_n\setminus B|}{n}
 \ge\frac{|\partial_e A_n|-|\partial_e B|}{4dn}
       -\frac{|B|-n}{2n}.
\end{equation}
On the probability one event in Theorem~\ref{thm:boundary}, we thus 
have that for every sequence of finite lattice sets $B_n$ with
$|B_n|=n+o(n)$ and $|\partial_e B_n|=o(n)$, 
\begin{equation}\label{eq:sublinear-boundary}
 \liminf_{n\to\infty}\frac{|A_n\setminus B_n|}{n}
 \ge\alpha\coloneqq\frac{b_d}{4d}>0
\end{equation}
as claimed in \eqref{eq:exclude-small-bdry} (where since the
comparison is deterministic, such sets $B_n$ may depend on
the entire \DLA{} realization).

By the preceding, for \eqref{eq:exclude-shape} it suffices to verify 
that $R_n$ of \eqref{def:Rn} has the required volume and boundary sizes. 
To this end, first fix a bounded Jordan measurable region $S$ of unit
volume, so $S \subset \R^d$ is bounded with $|\partial S|=0$ and let
$B_t\coloneqq\Z^d\cap tS$. Riemann sums give
\[
 |B_t|=t^d+o(t^d).
\]
The occupied endpoint of every edge in $\partial_e B_t$ is within
distance $1$ of $t\partial S$. Surround these endpoints by disjoint
unit lattice cells and rescale by $t^{-1}$. The cells lie within
distance $c_d/t$ of $\partial S$, so
\[
 |\partial_e B_t|
 \le 2d\,t^d
 \bigl|\{x\in\R^d\mid\dist(x,\partial S)\le c_d/t\}\bigr|
 =o(t^d)
\]
(where the last equality uses compactness and the measure zero of
$\partial S$). In particular, taking $S=R$ and $t=n^{1/d}$ establishes
\eqref{eq:exclude-shape} for all bounded regions of unit volume whose
boundaries have measure zero.

Now let $R$ be any compact set of unit volume. For any $\varepsilon>0$,
outer regularity and compactness give a finite union $D$ of bounded
open boxes with $R \subset D$ and $|D|<1+\varepsilon$.
Set $B_n\coloneqq\Z^d\cap (n^{1/d}D)$. Then
\[
 R_n\subset B_n,\qquad |B_n|=|D|n+o(n),\qquad
 |\partial_e B_n|=o(n).
\]
Applying \eqref{eq:shape-boundary-comparison} gives
\[
 \liminf_{n\to\infty}\frac{|A_n\setminus R_n|}{n}
 \ge\frac{b_d}{4d}-\frac{|D|-1}{2}
 \ge\alpha-\frac{\varepsilon}{2}.
\]
Letting $\varepsilon\downarrow0$ proves \eqref{eq:exclude-shape} for every
such compact $R$ of unit volume.
\qed

\section{Almost-double contacts: Proof of Theorem~\ref{thm:attachment}}\label{sec:attachment}

We start with the uniform attachment bound of part \textup{(a)} 
then we deduce the average-degree bound of part \textup{(b)}.
Fixing $0<\eta<1$, recall 
from \eqref{dfn:good-attachment} that 
\begin{align}\label{dfn:union-Ha}
E_\eta(C)&\,=\partial\mathcal O_C \cap\bigcup_{a\in C} H_a \\
H_a &\coloneqq
 \bigcup_{b\in C\setminus\{a\}}
 \{x\in\partial B(a,2)\colon 2 \le |x-b|\le 2+ \eta\}.
 \nonumber
\end{align}
Fixing $a \in C$, we proceed to uniformly bound the probability 
that $W(T_a) \in H_a$ for 
\[
 T_a \coloneqq\inf\{t\ge0\colon W(t)\in\overline B(a,2)\},
\]
by bounding the contribution of a single 
$b \in C$ followed by a union bound over the at most $6^d$
contributing neighbors, as illustrated in 
Figure~\ref{fig:Ha-two-steps}.

\begin{figure}
\centering
\begin{tikzpicture}[scale=.5,line cap=round,line join=round,
  dashed circle/.style={black,dashed,line width=.8pt}]
  \definecolor{thinred}{RGB}{180,38,45}
  \pgfmathsetmacro{\lo}{60}
  \pgfmathsetmacro{\hi}{acos((8-(2+.12)^2)/8)}
  \foreach \panel in {0,1}{
    \begin{scope}[xshift=\panel*10cm]
      \path[use as bounding box] (-4.2,-4.2) rectangle (4.2,4.2);
      \ifnum\panel=0
        \def\neighbors{10}
      \else
        \def\neighbors{10,95,190,280}
      \fi
      \filldraw[fill=gray!12,draw=black!70,line width=.65pt]
        (0,0) circle (1);
      \foreach \angle in \neighbors {
        \filldraw[fill=gray!12,draw=black!70,line width=.65pt]
          (\angle:2) circle (1);
      }
      \foreach \angle in \neighbors {
        \draw[dashed circle] (\angle:2) circle (2);
      }
      \draw[dashed circle] (0,0) circle (2);
      \foreach \angle in \neighbors {
        \foreach \sgn in {-1,1}{
          \pgfmathsetmacro{\start}{\angle+\sgn*\lo}
          \pgfmathsetmacro{\stop}{\angle+\sgn*\hi}
          \draw[thinred,line width=2.4pt] (\start:2)
            arc[start angle=\start,end angle=\stop,radius=2];
        }
      }
      \fill (0,0) circle (1.3pt) node[below left] {$a$};
      \ifnum\panel=0
        \fill (10:2) circle (1.3pt) node[below right] {$b$};
      \fi
    \end{scope}
  }
\end{tikzpicture}
\caption{
In the left panel, the red arcs show the contribution of one neighboring
center $b$ to $H_a$ i.e., the points $x\in\partial B(a,2)$
with $2\le |x-b|\le 2+\eta$. The right panel shows $H_a$ as the union
of these contributions over at most $6^d$ centers.
}
\label{fig:Ha-two-steps}
\end{figure}

\begin{lemma}\label{eu-lem:thin}
There are constants $c_d<\infty$ such that, for every
finite nonempty $2$-separated set $C$, every $a\in C$, and
$0<\eta<1$,
\[
 \sup_{|y-a|=3}\PP_y(T_a<\infty,\ W(T_a)\in H_a)
 \le c_d \, \sqrt\eta.
\]
When $d\ge3$, the \abbr{rhs} improves to $c_d \, \eta$.
\end{lemma}

\begin{proof} Write $\operatorname{SA}(\cdot)$ for surface area of a set, interpreted as arc length when $d=2$.
For each $b\in C\setminus\{a\}$, consider the spherical region
\begin{equation}\label{dfn:spherical-a-b}
 \{x\in\partial B(a,2)\colon 2 \le |x-b| \le 2+ \eta\}.
\end{equation}
Squaring the inequalities in \eqref{dfn:spherical-a-b}
shows that this region lies
between two parallel hyperplanes at distance at most $4 \eta/|a-b| \le 2\eta$
(since $|a-b|\ge2$). Such a strip on a sphere of
radius $2$ has surface area at most $c_d\sqrt\eta$ for $d=2$ and $c_d \eta$ for $d\geq 3.$

By the triangle inequality, the region in \eqref{dfn:spherical-a-b} is empty 
unless $|a-b| \le 4+\eta <5$. All open unit balls centered at such
values of $b$ lie inside $B(a,6)$ which fits at most
$6^d$ disjoint unit balls. Summing over their centers yields
\[
 \operatorname{SA}(H_a)\le c_{d} \, \sqrt\eta.
\]
When $d\ge3$, the \abbr{rhs} improves to $c_{d} \, \eta$.

By Poisson's formula \cite[Theorem~3.44]{Morters2010Brownian}, the density of the
unconditioned first hitting distribution on $\partial B(a,2)$
with respect to surface area is bounded, uniformly over
all starting points on $\partial B(a,3)$. Hence
\[
 \sup_{|y-a|=3}
 \PP_y(T_a<\infty,\ W(T_a)\in H_a)
 \le c_d\,\operatorname{SA}(H_a),
\]
which given our surface-area bounds completes the proof of the lemma.
\end{proof}

\begin{proof}[Proof of Theorem~\ref{thm:attachment}]$~$\\
(a). Fix a finite nonempty $2$-separated set $C$ and $0<\eta<1$.
Starting a Brownian motion at $z$ with $\dist(z,C)\ge3$, let
\[
 \sigma_a\coloneqq\inf\{t\ge0\colon |W(t)-a|\le3\}.
\]
Denote by $N$ the number of radius-$3$ neighborhoods around points in $C$ 
which the Brownian path enters before its first entrance $\tau_C$ to 
$\mathcal O_C=\bigcup_{a\in C}\overline B(a,2)$. That is, 
\begin{equation}\label{dfn:3-balls-count}
 N\coloneqq\sum_{a\in C}\ind_{\{\sigma_a<\tau_C\}},
 \qquad
 \EE_z [N] =\sum_{a\in C}\PP_z(\sigma_a<\tau_C).
\end{equation}
Since attachment in $H_a$ requires $\sigma_a<\tau_C=T_a$, 
the strong Markov property at $\sigma_a$ and Lemma~\ref{eu-lem:thin} yield the bound
\[
 \PP_z(\tau_C<\infty,\ W(\tau_C)\in H_a)
 \le c_d\sqrt\eta\,\PP_z(\sigma_a<\tau_C).
\]
In view of \eqref{dfn:union-Ha} and \eqref{dfn:3-balls-count}, 
summing the preceding bound over $a \in C$
yields that 
\begin{equation}\label{eq:bd-E-eta}
 \PP_z(\tau_C<\infty,\ W(\tau_C)\in E_\eta(C))
 \le c_d\sqrt\eta\,\EE_z[N] \,.
\end{equation}
Next, recall that by translation and rotation invariance, 
the probability $p_d>0$ that Brownian motion started on 
$\partial B(a,3)$ hits $\overline B(a,2)$ before leaving $B(a,4)$
depends only on $d$. Moreover, having $\sigma_a<\tau_C$ followed 
by the preceding event, implies that $\tau_C < \infty$ and
$W(\tau_C) \in B(a,4)$, regardless of where exactly 
the Brownian motion hits $\partial \mathcal O_C$. Thus, by 
the strong Markov property at $\sigma_a$, 
\begin{equation}\label{eq:bd-hit-3-ball-first}
 p_d\,\PP_z(\sigma_a<\tau_C)
 \le\PP_z(\tau_C<\infty,\ W(\tau_C)\in B(a,4)).
\end{equation}
Further, as $C$ is $2$-separated, any $x \in \R^d$ lies within distance $4$ 
of at most $5^d$ points $a \in C$ (by volume bounds, since the 
unit balls centered at such points of $C$ are disjoint and inside $B(x,5)$). Hence, 
\[
M\coloneqq \ind_{\{\tau_C<\infty\}} \sum_{a\in C}\ind_{\{W(\tau_C) \in B(a,4)\}} 
\le 5^d \ind_{\{\tau_C<\infty\}}
\]
and summing \eqref{eq:bd-hit-3-ball-first} over 
$a\in C$ yields, in view of \eqref{dfn:3-balls-count},
\[
 p_d\,\EE_z [N] \le\sum_{a\in C}
       \PP_z(\tau_C<\infty,\ W(\tau_C)\in B(a,4)) =\EE_z [M] \le 5^d\,  \PP_z(\tau_C<\infty).
\]
Combining the latter bound and \eqref{eq:bd-E-eta}
we arrive at 
\[
 \PP_z(\tau_C<\infty,\ W(\tau_C)\in E_\eta(C))
 \le\frac{5^d}{p_d} \, c_d \, \sqrt\eta\,
       \PP_z(\tau_C<\infty).
\]
Dividing both sides by $\PP_z(\tau_C<\infty)>0$ and taking $|z|\to\infty$
we arrive at the stated bound of \eqref{eq:unif-attachment}
(recall the definition \eqref{def:harmonic-measure} of $\mu_C(\cdot)$).

\noindent
(b). Fixing $\eta \in (0,1)$, consider the filtration $\cF_k \coloneqq\sigma(C_1,\ldots,C_k)$ and the 
adapted 
\[
\xi_k \coloneqq |\mathcal{E}(G_\eta(C_k))| - |\mathcal{E}(G_\eta(C_{k-1}))| 
\]
which counts for $k \ge 2$ the number of edges added to the proximity graph
$G_\eta(C_{k-1})$ by attaching $X_k$ to the Euclidean \DLA{}. With $2+\eta \le 3$,
it follows by volume estimates that any $2$-separated set has 
at most $4^d$ points in any ball of the form $\overline B(x,2+\eta)$,
so in particular $\xi_k \le 4^d$. Also, by definition $\xi_k = 1$
outside the event $\{X_k \in E_\eta(C_{k-1})\}$, hence 
we deduce from part (a) that 
\[
 \EE[\xi_k\mid\cF_{k-1}] \le 1+(4^d-1)\mu_{C_{k-1}}(E_\eta(C_{k-1}))
 \le 1+ (4^d-1) c_d \, \sqrt{\eta} \,.
\]
Consequently, 
\begin{equation}\label{eq:bd-av-degree}
\frac{|\mathcal{E}(G_\eta(C_n))|}{n} = \frac{1}{n} \sum_{k=2}^n \xi_k 
\le 1+ (4^d-1) c_d \, \sqrt{\eta} + n^{-1} S_n \,,
\end{equation}
for the $(\cF_n)$ martingale
\[
S_n \coloneqq \sum_{k=2}^n \bigl(\xi_k-\EE[\xi_k\mid\cF_{k-1}]\bigr) \,.
\]
Applying Azuma--Hoeffding to this martingale of uniformly bounded differences,
it follows by Borel--Cantelli that $n^{-1} S_n \to 0$ almost surely and upon 
multiplying $c_d$ by $4^d$ we get \eqref{eu-eq:degree-bound} 
out of \eqref{eq:bd-av-degree}.

Lastly, for $d\ge 3$ the sharper estimate in Lemma~\ref{eu-lem:thin}
allows us to replace $\sqrt\eta$ by $\eta$ throughout the proof of both 
parts of the theorem.
\end{proof}

\section{Excluding Optimal Packing: Proof of Corollary~\ref{cor:euclidean}}\label{sec:euclidean}

Recall that a periodic packing is an infinite $2$-separated set $\mathcal P\subset\R^d$ 
obtained from $m=|\mathcal P\cap Q| \ge 1$ points inside a fundamental box
\begin{equation}\label{def:fund-box}
 Q=\prod_{i=1}^d[0,T_i),
\end{equation}
which are then repeated with periods $T_1,\ldots,T_d$ along the coordinate axes. The number density $r(\mathcal P)$ and the average proximity degree
$\overline d_\eta(\mathcal P)$ are given by
\begin{equation}\label{dfn:d-eta-P}
 r(\mathcal P)=\frac{m}{|Q|},\qquad
 \overline d_\eta(\mathcal P)=\frac1m
 \sum_{v\in\mathcal P\cap Q} \deg_{G_\eta(\mathcal P)} (v)
\end{equation}
Our next lemma shows that in the
plane a periodic packing with average proximity degree at most five
has number density bounded away from the optimal number density
$\rho_d\coloneqq\Delta_d/\omega_d$.
\begin{lemma}\label{eu-lem:localpacking}
Let $d=2$ and $\eta \in (0,1/4)$. There exists 
$\varepsilon=\varepsilon(d,\eta)>0$ such that for every 
periodic packing in $\R^d$,
\begin{equation}\label{eu-eq:periodic-gap}
 \overline d_\eta(\mathcal P)\le 5
 \quad\Longrightarrow\quad
 r(\mathcal P)\le (1-2\varepsilon) \rho_d.
\end{equation}
\end{lemma}

\begin{proof}
\hypertarget{eu-proof:localpacking}{}
The Voronoi cells of points $v$ in a periodic packing $\mathcal P\subset\R^2$,
\[
 \Omega_v \coloneqq \{x\in\R^2\colon |x-v|\le |x-w|
                          \text{ for every }w\in\mathcal P\},
\]
are bounded convex polygons containing the corresponding unit disks.
They form a partition of the plane up to a boundary set of zero area.
Recall from \cite{Mazel2022Minimal} that
$|\Omega_v|\ge 2\sqrt3=\rho_2^{-1}$ for every $v\in\mathcal P$,
with equality only if $\Omega_v$ is a regular hexagon circumscribed
about the unit disk centered at $v$. In this case, $v$ has six
neighboring centers at distance $2$.

Since the periodic extensions by $\{T_i\}$ of both $Q$ and 
$\{\Omega_v,  v \in \mathcal P \cap Q\}$ partition $\R^2$, it follows 
by Fubini's theorem that 
\begin{equation}\label{eq:vol-id}
|Q| = \sum_{v \in \mathcal P \cap Q} |\Omega_v| \,.
\end{equation}
 Having $\overline d_\eta(\mathcal{P}) \le 5$ implies, see \eqref{dfn:d-eta-P},
that $\deg_{G_\eta(\mathcal P)} (v)
\le 5$ for at least $m/6$ points $v \in \mathcal P \cap Q$,
so in view of \eqref{eq:vol-id} it suffices to show that for 
some $a_\eta>0$ and all $v \in \mathcal P$,
\begin{equation}\label{eu-eq:local-excess}
\deg_{G_\eta(\mathcal P)} (v)
 \le5
 \quad\Longrightarrow\quad
 |\Omega_v|\ge2\sqrt3+a_\eta \,.
\end{equation}

As illustrated in Figure~\ref{fig:planar-compactness-three} this 
holds because the center $v$ of 
any sequence  of cells with $|\Omega_v| \to 2 \sqrt{3}$ must acquire 
along the way six near neighbors.
\begin{figure}[t]
\centering
\begin{tikzpicture}[scale=.61,line cap=round,line join=round]
\definecolor{cellblue}{RGB}{45,99,140}
\definecolor{excessred}{RGB}{173,75,48}
\begin{scope}[xshift=0.0cm]
\draw[black!40,line width=.45pt] (60:2) circle (1);
\draw[black!40,line width=.45pt] (120:2) circle (1);
\draw[black!40,line width=.45pt] (180:2) circle (1);
\draw[black!40,line width=.45pt] (240:2) circle (1);
\draw[black!40,line width=.45pt] (300:2) circle (1);
\draw[black!40,line width=.65pt] (2,0) circle (1);
\fill[cellblue!17] (1.000000,0.577350) -- (0.000000,1.154701) -- (-1.000000,0.577350) -- (-1.000000,-0.577350) -- (-0.000000,-1.154701) -- (1.000000,-0.577350) -- cycle;
\draw[cellblue,line width=.85pt] (1.000000,0.577350) -- (0.000000,1.154701) -- (-1.000000,0.577350) -- (-1.000000,-0.577350) -- (-0.000000,-1.154701) -- (1.000000,-0.577350) -- cycle;
\draw[black,line width=.6pt] (0,0) circle (1);
\draw[black!65,dashed,line width=.55pt] (0,0) circle (2.2);
\node[fill=white,inner sep=1pt,font=\small] at (90:2.2) {$2+\eta$};
\fill ({60}:2) circle (1.7pt);
\fill ({120}:2) circle (1.7pt);
\fill ({180}:2) circle (1.7pt);
\fill ({240}:2) circle (1.7pt);
\fill ({300}:2) circle (1.7pt);
\fill[black] (2,0) circle (2pt);
\fill (0,0) circle (1.7pt);
\node[below left,inner sep=1pt,font=\small] at (0,0) {$v$};
\node[font=\small] at (0,-3.45) {(a) $|\Omega_v|=2\sqrt3$};
\end{scope}
\begin{scope}[xshift=7.2cm]
\draw[black!40,line width=.45pt] (60:2) circle (1);
\draw[black!40,line width=.45pt] (120:2) circle (1);
\draw[black!40,line width=.45pt] (180:2) circle (1);
\draw[black!40,line width=.45pt] (240:2) circle (1);
\draw[black!40,line width=.45pt] (300:2) circle (1);
\draw[excessred,line width=.65pt] (2.65,0) circle (1);
\fill[cellblue!17] (1.325000,0.389711) -- (0.000000,1.154701) -- (-1.000000,0.577350) -- (-1.000000,-0.577350) -- (-0.000000,-1.154701) -- (1.325000,-0.389711) -- cycle;
\fill[excessred!45] (1,0.5773502691896258) -- (1.325,0.38971143170299744) -- (1.325,-0.38971143170299744) -- (1,-0.5773502691896258) -- cycle;
\draw[black!45,densely dotted] (1,0.5773502691896258) -- (1,-0.5773502691896258);
\draw[cellblue,line width=.85pt] (1.325000,0.389711) -- (0.000000,1.154701) -- (-1.000000,0.577350) -- (-1.000000,-0.577350) -- (-0.000000,-1.154701) -- (1.325000,-0.389711) -- cycle;
\draw[black,line width=.6pt] (0,0) circle (1);
\draw[black!65,dashed,line width=.55pt] (0,0) circle (2.2);
\node[fill=white,inner sep=1pt,font=\small] at (90:2.2) {$2+\eta$};
\fill ({60}:2) circle (1.7pt);
\fill ({120}:2) circle (1.7pt);
\fill ({180}:2) circle (1.7pt);
\fill ({240}:2) circle (1.7pt);
\fill ({300}:2) circle (1.7pt);
\fill[excessred] (2.65,0) circle (2pt);
\fill (0,0) circle (1.7pt);
\node[below left,inner sep=1pt,font=\small] at (0,0) {$v$};
\node[font=\small] at (0,-3.45) {(b) $|\Omega_v|\ge2\sqrt3+a_\eta$};
\end{scope}
\begin{scope}[xshift=14.4cm]
\draw[black!40,line width=.45pt] (60:2) circle (1);
\draw[black!40,line width=.45pt] (120:2) circle (1);
\draw[black!40,line width=.45pt] (180:2) circle (1);
\draw[black!40,line width=.45pt] (240:2) circle (1);
\draw[black!40,line width=.45pt] (300:2) circle (1);
\draw[excessred,line width=.65pt] (2.12,0) circle (1);
\fill[cellblue!17] (1.060000,0.542709) -- (0.000000,1.154701) -- (-1.000000,0.577350) -- (-1.000000,-0.577350) -- (-0.000000,-1.154701) -- (1.060000,-0.542709) -- cycle;
\fill[excessred!45] (1,0.5773502691896258) -- (1.06,0.5427092530382482) -- (1.06,-0.5427092530382482) -- (1,-0.5773502691896258) -- cycle;
\draw[black!45,densely dotted] (1,0.5773502691896258) -- (1,-0.5773502691896258);
\draw[cellblue,line width=.85pt] (1.060000,0.542709) -- (0.000000,1.154701) -- (-1.000000,0.577350) -- (-1.000000,-0.577350) -- (-0.000000,-1.154701) -- (1.060000,-0.542709) -- cycle;
\draw[black,line width=.6pt] (0,0) circle (1);
\draw[black!65,dashed,line width=.55pt] (0,0) circle (2.2);
\node[fill=white,inner sep=1pt,font=\small] at (90:2.2) {$2+\eta$};
\fill ({60}:2) circle (1.7pt);
\fill ({120}:2) circle (1.7pt);
\fill ({180}:2) circle (1.7pt);
\fill ({240}:2) circle (1.7pt);
\fill ({300}:2) circle (1.7pt);
\fill[excessred] (2.12,0) circle (2pt);
\fill (0,0) circle (1.7pt);
\node[below left,inner sep=1pt,font=\small] at (0,0) {$v$};
\draw[->,excessred,line width=.7pt] (2.85,1.3) -- (2.05,1.3);
\node[font=\small] at (0,-3.45) {(c) $|\Omega_v|\longrightarrow2\sqrt3$};
\end{scope}
\end{tikzpicture}
\caption{The shaded polygons are Voronoi cells of the central unit disks.
(a) Six contacts give the minimum area $2\sqrt3$.
(b) Moving the sixth center beyond the dashed circle of radius $2+\eta$
creates an area excess, shown in orange.
(c) As the excess vanishes, the sixth disk center must enter the dashed circle before contact.}
\label{fig:planar-compactness-three}
\end{figure}
Indeed, if \eqref{eu-eq:local-excess} fails, there is a sequence of
$2$-separated periodic packings $(\mathcal P_n)_{n\ge1}$, each with
a chosen center having at most five near neighbors, such that the
areas of the corresponding Voronoi cells tend to $2\sqrt3$.
Translate each packing so that its chosen center lies at $0$, and
denote the cell at $0$ by
\begin{equation}\label{dfn:Omega-n}
 \Omega^{(n)}\coloneqq
 \{x\in\R^2\colon |x|\le |x-w|\text{ for every }w\in\mathcal P_n\}.
\end{equation}
By construction
\begin{equation}\label{eq:packing-contradiction}
 \bigl|(\mathcal P_n\setminus\{0\})\cap\overline B(0,2+\eta)\bigr|
 \le 5
 \qquad\text{and}\qquad
 |\Omega^{(n)}|\to 2\sqrt3 \quad\text{as }n\to\infty.
\end{equation}

Since $\overline B(0,1) \subset \Omega^{(n)}$, 
for any unit vector $u$ perpendicular to some $x\in\Omega^{(n)}$ 
the convex set $\Omega^{(n)}$ must contain the triangle 
with vertices $x,u,-u$ whose area is $|x|$. Consequently, 
$\Omega^{(n)} \subset \overline B(0,4)$ once $|\Omega^{(n)}| \le 4$,
whereupon it suffices by the triangle inequality to consider in 
\eqref{dfn:Omega-n} only $w \in \mathcal P_n \cap \overline B(0,8)$. 
Further, due to volume considerations there are at most 
$9^2$ such points $w$ in the $2$-separated $\mathcal P_n$.
By compactness, passing to a subsequence $n_k$, the $2$-separated collections 
$\mathcal P_{n_k} \cap \overline B(0,8)$ converge, as do the 
corresponding Voronoi cells $\Omega^{(n_k)}$ they fully determine. 
The limit of $\Omega^{(n_k)}$ must be a Voronoi cell 
in the limiting configuration, of area $2 \sqrt{3}$ and center $0$.  
As such it must be a regular hexagon with $0$ having six neighbors 
at distance $2$ in the limiting configuration. The latter 
correspond in turn to having $6$ points in 
$\bar B(0,2+\eta) \cap \mathcal P_{n_k} \setminus \{0\}$
for all sufficiently large $k$, contradicting \eqref{eq:packing-contradiction}.
\end{proof}

\begin{remark}\label{rmk:sphere-packing}
We believe that Lemma~\ref{eu-lem:localpacking}, and hence
Corollary~\ref{cor:euclidean}, also holds in dimensions
$d\in\{3,8,24\}$, where optimal sphere packings are well understood
\cite{Hales2005Kepler,Borroczky2024Stability}. 
The arguments in these dimensions appear to require
substantial machinery from sphere packing theory. We thus  
restrict to the planar case in order to keep the 
focus on our probabilistic results.
\end{remark}

\begin{proof}[Proof of Corollary~\ref{cor:euclidean}]
We give the argument keeping $d\ge2$ general since it works so long as Lemma~\ref{eu-lem:localpacking} holds in that dimension. We then deduce our corollary by setting $d=2$.
For $\eta>0$ and a finite, non-empty $2$-separated set $Y$, write
\begin{align}\label{dfn:d-eta}
\overline d_\eta(Y) &\coloneqq \frac{2}{|Y|}|\mathcal E(G_\eta(Y))|, \qquad 
\end{align}
for the average degree of the associated proximity graph. Fix a
compact set $R \subset \R^d$ of unit volume and set 
\begin{equation}\label{dfn:rho-eta}
\begin{split}
\rho_{\eta,R}(t)\coloneqq t^{-d}\sup\bigl\{|Y|\colon {}& |Y|\ge1,\ Y\subset tR,\\
 &Y\text{ is a $2$-separated set},\ \overline d_\eta(Y)\le5\bigr\}.
\end{split}
\end{equation}
Assuming Lemma~\ref{eu-lem:localpacking} holds in dimension $d$,
we will deduce from \eqref{eu-eq:periodic-gap} that for every
$\eta\in(0,1/4)$
\begin{equation}\label{eq:rho-sub-optimal}
\limsup_{t \to \infty} \rho_{\eta,R} (t) \le (1-\varepsilon) \rho_d.
\end{equation}
Indeed, partitioning $\R^d$ into half-open cubes of side $L$ such that 
\[
 (1-2 \varepsilon)(1+3/L)^d \le 1-\varepsilon \,,
\]
the union of the $J_t$ cubes meeting $t R$ contains $tR$ and lies within distance $\sqrt d L$ of it.
Hence, by volume rescaling and applying downward continuity of Lebesgue measure to nested open neighborhoods of the compact $R$,
\begin{equation}\label{eu-eq:grid-volume}
\lim_{t \to \infty} t^{-d} J_t L^d= 1.
\end{equation}
We shall associate with any $2$-separated $Y\subset tR$
a periodic packing $\mathcal P$ so that 
\begin{equation}\label{eq:r-P}
 r(\mathcal P)=\frac{|Y|}{J_t(L+3)^d},\qquad
 \overline d_\eta(\mathcal P)\le\overline d_\eta(Y).
\end{equation}
This would imply by \eqref{eu-eq:periodic-gap}, \eqref{eu-eq:grid-volume} and the preceding choice of $L$, that 
\[
\limsup_{t \to \infty} \rho_{\eta,R} (t) 
\le (1-2\varepsilon) \rho_d (1+3/L)^d \limsup_{t \to \infty} \{t^{-d} J_t L^d\}
 \le(1-\varepsilon)\rho_d\,,
\]
as claimed in \eqref{eq:rho-sub-optimal}. To get \eqref{eq:r-P}, note
that to periods $T_1=J_t (L+3)$, $T_i=(L+3)$, $i \ge 2$,
corresponds a box of volume $|Q|=J_t (L+3)^d$. Recalling our $J_t$ disjoint 
translates of $[0,L)^d$ by $\{\vec v_i \}$,
set the $m=|Y|$ points in $\mathcal P \cap Q$ as follows:
\[
\mathcal P \cap Q \coloneqq
\bigcup_{i=1}^{J_t}
\bigl\{ (i-1)(L+3)e_1 + y \colon
        y \in (Y-\vec v_i)\cap[0,L)^d \bigr\}
\]
(that is, we retain the positions of points of $Y$ relative to the lower corner 
of the side $L$ cube which they intersect, putting those $J_t$ patterns 
in the $[0,L)^d$-sub-cubes within a row of $J_t$ boxes of side $L+3$ each). The
number density of \eqref{eq:r-P} is obvious, with the 
\abbr{rhs} of \eqref{eq:r-P} holding because the gaps of width $3$ 
prevent forming any new pair touching $\mathcal P \cap Q$ and having 
distance at most $2+\eta$.

Choose $\eta \in (0,1/4)$ so the \abbr{rhs} of \eqref{eu-eq:degree-bound} is at most $3/2$. 
By Theorem~\ref{thm:attachment}\textup{(b)}, 
excluding a single event of probability zero, we can fix 
a realization of the Euclidean \DLA{} such that  
$\overline d_\eta(C_n) \le 4$ for some $n_0$ finite and all $n \ge n_0$. 
Let $R \subset \R^d$ be a compact set of unit volume, possibly 
dependent on the realization. Setting $Y_n \coloneqq C_n \cap P_n(R)$
recall \eqref{dfn:PnR} that $P_n(R)=t_n R$ for $t_n^d=n/\rho_d$. 
Then, with $|C_n|=n$, for $n \ge n_0$, by definitions \eqref{dfn:d-eta} and \eqref{dfn:rho-eta}
\[
\frac{|Y_n|}{n} \ge \frac{4}{5} \quad \Longrightarrow \quad 
\overline d_\eta(Y_n) \le \frac{|C_n|}{|Y_n|} \overline d_\eta (C_n) \le 5 
\quad \Longrightarrow \quad \frac{|Y_n|}{n} \le \frac{\rho_{\eta,R}(t_n)}{\rho_d} \,.
\]
We thus deduce from \eqref{eq:rho-sub-optimal} that 
\[
\limsup_{n \to \infty} \frac{|Y_n|}{n} \le \max\{\frac{4}{5},1-\varepsilon\},
\]
which proves the corollary for $\alpha=\min\{\frac{1}{5},\varepsilon\}$.
\end{proof}

\subsection*{Acknowledgments}

Dembo was partially supported by NSF grant DMS-2348142. 
Junge was partially supported by NSF grant DMS-2115936.

\subsection*{Statement on AI use}
We began this project by asking ChatGPT to prove Corollary 1.2 in the special case of a ball and to prove the analogous result for Euclidean DLA. ChatGPT provided proof drafts of Corollaries 1.2 and 1.5 for this special case as well as a proof of Theorem 1.4. Once we understood the key ideas, we explored possible extensions with ChatGPT, which led to Theorems 1.1 and 1.3 and the broader exclusion statements in Corollaries 1.2 and 1.5. 

We rewrote the manuscript ourselves, occasionally incorporating suggestions from ChatGPT. In doing so, we simplified the description of $K$ and wrote the energy martingale arguments in a more standard form. For Euclidean DLA, we devised a more direct proof of Theorem 1.4 using the Poisson kernel and simplified the packing argument in the proof of Corollary 1.5 to make it more accessible to probabilists. In addition to ironing out many smaller details, we deduced shape exclusion from linear boundary growth and observed that the estimates for permanent vacancies yield a local boundary bound in cubes. 

\bibliographystyle{alpha}
\bibliography{DLA}
\end{document}